\documentclass[12pt]{article}
\usepackage{color}
\usepackage[fleqn]{amsmath}
\usepackage{amscd}
\usepackage{amsfonts}
\usepackage{amssymb}
\usepackage{amsthm}
\usepackage{bm}
\usepackage{float}
\usepackage{graphicx}
\usepackage{epstopdf}
\renewcommand{\figurename}{Fig.} 
\usepackage{subfigure}
\usepackage{multirow} 
\usepackage{tabularx}
\usepackage{makecell} 
\usepackage[font=footnotesize]{caption} 

\usepackage{cite}

\allowdisplaybreaks 

\makeatletter
\newenvironment{breakablealgorithm}
{
		\begin{center}
			\refstepcounter{algorithm}
			\hrule height.8pt depth0pt \kern2pt
			\renewcommand{\caption}[2][\relax]{
				{\raggedright\textbf{\ALG@name~\thealgorithm} ##2\par}%
				\ifx\relax##1\relax 
				\addcontentsline{loa}{algorithm}{\protect\numberline{\thealgorithm}##2}%
				\else 
				\addcontentsline{loa}{algorithm}{\protect\numberline{\thealgorithm}##1}%
				\fi
				\kern2pt\hrule\kern2pt
			}
		}{
		\kern2pt\hrule\relax
	\end{center}
}
\makeatother

\usepackage[english]{babel}
\usepackage[utf8]{inputenc}
\usepackage{algorithm}
\usepackage{algorithmicx}

\usepackage[noend]{algpseudocode}

\usepackage[framemethod=tikz]{mdframed} 
\usepackage{lipsum}
\usepackage{soul, color, xcolor}

\usepackage{mathrsfs}

\newtheorem{theorem}{Theorem}[section]
\newtheorem{lemma}[theorem]{Lemma}

\newtheorem{remark}[theorem]{Remark}

\usepackage{indentfirst} 

\begin{document}
\baselineskip 20pt

\title{\bf Greedy randomized block Kaczmarz method
for matrix equation $AXB=C$
and its applications in color image restoration\
\thanks{The work is supported by the Fundamental Research Funds for the Central Universities (YA24YJS00040) and National Natural Science Foundation of China (No.61931003). 
}}

\author{\bf
Wenli Wang$^{a}${\thanks{E-mail addresses: 22110483@bjtu.edu.cn}}, Duo Liu$^{a}${\thanks{E-mail addresses: 20118001@bjtu.edu.cn}}, Gangrong Qu$^{a}${\thanks{E-mail addresses: grqu@bjtu.edu.cn}}
\\
\normalsize\sl $^{a}${\small {\it School of Mathematics and Statistics, Beijing Jiaotong University, Beijing 100044, China}}\\
}

\date{}
\maketitle

\begin{center}
\textbf{This is an old version. Please refer to the updated version available at https://arxiv.org/abs/2602.03239.
}
\end{center}

{\noindent}\\
------------------------------------------------------------
-----------------------------------------------------------\\

{\small \noindent
{\bf Abstract}\\

%

In view of the advantages of simplicity and effectiveness of the Kaczmarz method, which was originally employed to solve the large-scale system of linear equations $Ax=b$, we study the greedy randomized block Kaczmarz method (ME-GRBK) and its relaxation and deterministic versions to solve the matrix equation $AXB=C$, which is commonly encountered in the applications of engineering sciences. It is demonstrated that our algorithms converge to the unique least-norm solution of the matrix equation when it is consistent and their convergence rate is faster than that of the randomized block Kaczmarz method (ME-RBK). Moreover, the block Kaczmarz method (ME-BK) for solving the matrix equation $AXB=C$ is investigated and it is found that the ME-BK method converges to the solution $A^{+}CB^{+}+X^{0}-A^{+}AX^{0}BB^{+}$ when it is consistent. The numerical tests verify the theoretical results and
the methods presented in this paper are applied to the color image restoration problem to obtain satisfactory restored images.
\\

\noindent {\it Key words:}~~Greedy randomized block Kaczmarz method;
Maximal weighted residual block Kaczmarz method;
Color image restoration\\
\noindent {\it AMS subject classification:} 15A24; 15A21; 65F10

\section{Introduction}

Big data and the growing scale of engineering sciences problems have brought new challenges to the design of algorithms for large-scale matrix equations [1-4].
Direct methods for computing small or medium scale system of equations are frequently impractical for large-scale system of equations. As a result, iterative methods that seek an optimal solution using the objective function are gaining popularity.

Given the matrices $A\in \mathbb{R}^{m\times p}$, $B\in \mathbb{R}^{q\times n}$ and $C\in \mathbb{R}^{m\times n}$, where $m,n,p$ and $q$ may be large,
we consider to solve the matrix equation
\begin{align*}
AXB=C,
\tag{$1.1$}
\end{align*}
which has a close connection to numerous important applications in the engineering field, such as image processing, stability analysis, control theory and the regression analysis in machine learning [2-7].
Until now, a variety of high-performance methods have been offered for solving matrix equation (1.1) [3, 4, 8-14, 32].
However, large-scale computing problem may not be appropriate to be solved by some traditional centralised algorithms, such as gradient-based iterative (GI) algorithm [8], iterative orthogonal direction (IOD) method [9], Jacobi
and Gauss-Seidel-type
iteration methods [10] and other iterative methods that require all data to be stored and computed at each iteration. This puts a high demand on the storage and computing power of the computer.

The row-action and column-action methods have important feature that they do not
require access to the entire coefficient matrix at each iteration, avoiding the concern that the matrix-vector product is quite expensive or the coefficient matrix is too large to be stored, and are suitable for solving large-scale computing problems [15-29].
Recently, Zeng et al. [3, 4] studied the distributed computation over multiagent networks,
where each agent only accesses the subblock matrices of $A$, $B$ and $F$.
Specifically, they obtained the distributed structures using row-action and column-action so that the original problem (1.1) becomes some distributed optimisation problems with different constraints.
In addition,
coordinate descent algorithm is a typical column-action method and has been generalised to the randomized block coordinate descent (RBCD) method to solve the matrix equation (1.1) [13].


The classical Kaczmarz method [15]
\begin{align*}
x^{k+1}=x^{k}+\frac{b_{i_{k}}-A_{i_{k},:}x^{k}}{\|A_{i_{k},:}\|_{2}^{2}}
A_{i_{k},:}^{T},
\quad i_{k}=(k\ mod\ m)+1,
\tag{$1.2$}
\end{align*}
where $A_{i_{k},:}$ is the $i_{k}$-th row of $A$ and $b_{i_{k}}$ is the $i_{k}$-th component of $b$,
is a row-action method. It
was originally used to solve the large-scale linear system $Ax=b$.
In the field of medical imaging, Kaczmarz method is known as algebraic reconstruction technique (ART), which is the first iteration scheme based on block-row operation in Computed Tomography (CT) image reconstruction [16].

Due to its simplicity and efficiency, numerous improvement strategies for the Kaczmarz method have been proposed.
Strohmer and Vershynin [18] adopted a random probability criterion instead of a sequential manner for selecting a working row to participate in an iteration operation, resulting in a randomized Kaczmarz (RK) method with expected exponential rate of convergence, which attracted a great deal of attention and was further generalized to underdetermined or rank-deficient linear systems [21-24].
Subsequently, Bai and Wu [19] found that if the coefficient matrix $A$ is scaled by a diagonal matrix, then the norms of every row in $A$ equal to an identical constant, and the random probability criterion becomes uniform sampling.
To conquer this, they presented a greedy randomized Kaczmarz (GRK) method by introducing a more effective probability criterion.
Furthermore, the GRK method is extended by inserting a relaxation factor in probability criterion to obtain a relaxed greedy randomized Kaczmarz (RGRK) method [20]. Besides, other prominent improvements of RK method include the maximal weighted residual Kaczmarz (MWRK) method and the randomized extended Kaczmarz (REK) method; see, e.g., [25-28] and the references therein.

More lately, the extensions of Kaczmarz-type methods for solving the matrix equation (1.1)
has receive increased attention.
Based on the Petrov-Galerkin conditions, Wu et al. [11] investigated
a relaxed greedy randomized Kaczmarz (ME-RGRK) method and a maximal weighted residual Kaczmarz (ME-MWRK) method,
and showed that the presented methods have exponential convergence rates.
Soon after, a randomized block Kaczmarz (ME-RBK) method was developed in [12].
See also [30, 31]
for the generalization of Kaczmarz method to Sylvester matrix equation.



Based on [12] and inspired by [11, 19, 20, 25], we present a greedy randomized block Kaczmarz method (ME-GRBK), a relaxed greedy randomized block Kaczmarz method (ME-RGRBK), and a maximal weighted residual block Kaczmarz method (ME-MWRBK) for solving the large-scale matrix equation (1.1).
It is specifically stated that the algorithms proposed in this paper are different from the algorithms presented in [11]. In [11], the original matrix equation (1.1) transforms into $mn$ subsystems $A_{i,:}XB_{:,j}=C_{i,j}$ after row-action and column-action, and two indexes $i_{k},j_{k}$ need to be selected according to the stochastic probability criterion at each iteration. However, in this paper, the original matrix equation (1.1) becomes $m$ subsystems $A_{i,:}XB=C_{i,:}$ after row-action, and one index $i_{k}$ is picked via the stochastic probability criterion at each iteration. The main contribution in this paper are as follows.
\begin{itemize}
\item [$\bullet$]The convergence of the block Kaczmarz method (ME-BK) for solving the matrix equation (1.1) has not yet been investigated. This paper therefore sets out to explore this method and to demonstrate that the ME-BK method converges to the solution $A^{+}CB^{+}+X^{0}-A^{+}AX^{0}BB^{+}$ when the matrix equation (1.1) is consistent.
\item [$\bullet$]The ME-GRBK method is a greedy version of the ME-RBK method presented in [12], and the ME-RGRBK and ME-MWRBK methods can be regarded as the relaxation and deterministic versions of the ME-GRBK method.
    We prove that the ME-GRBK, ME-RGRBK and ME-MWRBK methods  converge to the unique least-norm solution $A^{+}CB^{+}$ when (1.1) is consistent and exhibit superior performance.
\item [$\bullet$]The numerical examples verify the theoretical results and illustrate that the ME-GRBK, ME-RGRBK and ME-MWRBK methods are more efficient than the ME-RBK method presented in [12].
    Furthermore, the methods have been applied to the color image restoration problem, resulting in satisfactory restored results.
\end{itemize}


The paper is organized as follows. Section 2 gives some necessary notations and preliminaries. Section 3 shows the ME-BK method and its convergence analysis. Section 4 investigates the ME-GRBK method and its relaxation and deterministic versions and their convergence properties. Section 5 reports some numerical examples to demonstrate the theoretical findings. Section 6 provides the application of our methods in color image restoration. Section 7 gives some concise conclusions.

\section{Notations and preliminaries}

For matrix $X\in \mathbb{R}^{m\times n}$, we use $X^{+}$, $X^{T}$, $r(X)$, $R(X)$, $\|X\|_{2}$, $\|X\|_{F}$, $X_{i,:}$, $X_{:,i}$ and $tr(X)$ 
to represent its Moore-Penrose inverse, the transpose, the rank, the range space, the spectral norm, the Frobenious norm, the $i$-th row, the $i$-th coloumn and the trace of $X$, respectively.
The inner product and the Kronecker product of matrices $A$ and $B$ are denoted as $\langle A,B \rangle$ and $A\otimes B$ respectively.
The stretching operator $vec(\cdot)$ is represented by $vec(X)=(x_{1}^{T},x_{2}^{T},\cdots,x_{n}^{T})^{T}\in \mathbb{R}^{mn}$ for matrix $X=(x_{1},x_{2},\cdots,x_{n})\in \mathcal{R}^{m\times n}$.
The identity matrix is denoted as $I$ and its size is known by the context. For two arbitrary integers $a$ and $b$, $(a\ mod\ b)$ denotes the remainder produced by dividing $a$ by $b$.
For arbitrary positive integer $m$, let $[m]=\{1,2,\cdots,m\}$.
The expected value conditional for the first $k$ iterations is denoted as $E_{k}$, i.e.,
\begin{align*}
E_{k}[\cdot]=E[\cdot|i_{0},i_{1},\cdots,i_{k-1}],
\end{align*}
where $i_{m}(m=0,1,\cdots,k-1)$ is the row index selected at the $m$-th iteration.

Let $\mathcal{A}=B^{T}\otimes A$, $x=vec(X)$ and $c=vec(C)$.
When (1.1) is inconsistent, we have $x\notin R(\mathcal{A})$
and are interested in calculating its least-squares solution
\begin{align*}
X_{LS}=\mathop{\arg\min}_{X\in\mathbb{R}^{p\times q}}\|C-AXB\|_{F}^{2}.
\tag{$2.1$}
\end{align*}
When (1.1) is consistent, it satisfies $x\in R(\mathcal{A})$
and has a unique solution or infinitely many solutions.
To be specific, if $r(\mathcal{A})= r(\mathcal{A},c)=pq$, then $A$ has full column rank, $B$ has full row rank, and (1.1) has a unique solution $(A^{T}A)^{-1}A^{T}CB^{T}(BB^{T})^{-1}$; if $r(\mathcal{A})= r(\mathcal{A},c)<pq$, then (1.1) has infinitely many solutions
\begin{align*}
X_{M}=A^{+}CB^{+}+Y-A^{+}AYBB^{+},\ \forall Y\in\mathbb{R}^{p\times q},
\tag{$2.2$}
\end{align*}
where $(\mathcal{A}, c)$ is an augmented matrix.
It is obvious that $A(Y-A^{+}AYBB^{+})B=\mathbf{0}$.
Let $X_{\ast}:=A^{+}CB^{+}$. Then $X_{\ast}$ is the unique least-norm solution, i.e., $X_{\ast}=\mathop{\arg\min}\limits_{AXB=C}\|X\|_{F}$.
This is because $\left\langle X_{\ast}, Y-A^{+}AYBB^{+}\right\rangle=0$ and
\begin{align*}
\|X_{M}\|_{F}^{2}
&=\|X_{\ast}\|_{F}^{2}+\|Y-A^{+}AYBB^{+}\|_{F}^{2}
\geq\|X_{\ast}\|_{F}^{2}.
\end{align*}
Specially, if $r(\mathcal{A})= r(\mathcal{A},c)=mn< pq$, then $A$ has full row rank, $B$ has full column rank, and $A^{T}(AA^{T})^{-1}C(B^{T}B)^{-1}B^{T}$ is the unique least-norm solution.

We give the following valuable lemma which will be utilized in the subsequent sections, and its proof can be found in [11, 31].

\begin{lemma} 
Let $A$ be an arbitrary non-zero real matrix. For any $x\in R(A)$, it holds
\begin{align*}
\|A^{T}x\|_{2}^{2}\geq \sigma_{\min}^{2}(A)\|x\|_{2}^{2}
=\frac{\|A\|_{F}^{2}}{\kappa_{F,2}^{2}(A)}\|x\|_{2}^{2},
\end{align*}
where the condition number is defined as $\kappa_{_{F,2}}(A):=\|A\|_{F}\|A^{+}\|_{2}=\|A\|_{F}/\sigma_{\min}(A)$.

\end{lemma}

The following ME-RBK method was proposed by Xing et al. in [12]. It was proved in [12] that the ME-RBK method is convergent in expectation to the least-norm solution $X_{\ast}$ when the matrix equation (1.1) is consistent, as in Lemma 2.2 below.
\begin{breakablealgorithm}
\caption{~ME-RBK method}
{ \begin{algorithmic}[1]
\noindent {\bf Initialization:} Give any initial iterative matrix $X^{0}$.
\State {\bf for} $k=0,1,2,\cdots$ until convergence {\bf do}
\State ~~~~Select $i_{k}\in\{1,2,\cdots,m\}$ with probability
${\rm Pr}({\rm row=}i_{k})
=\frac{\|A_{i_{k},:}\|_{2}^{2}}{\|A\|_{F}^{2}}$
\State ~~~~Compute $X^{k+1}=X^{k}+\frac{\alpha}{\|A_{i_{k},:}\|_{2}^{2}}
A_{i_{k},:}^{T}\left(
C_{i_{k},:}-A_{i_{k},:}X^{k}B\right)B^{T}$
\State {\bf end for}
\end{algorithmic}}
\end{breakablealgorithm}

\begin{lemma}
Let the matrix equation (1.1) with $A\in \mathbb{R}^{m\times p}$, $B\in \mathbb{R}^{n\times q}$ and $C\in \mathbb{R}^{m\times n}$
be consistent.
Assume that $0<\alpha<\frac{2}{\|B\|_{2}^{2}}$,
then the sequence $\{X^{k}\}_{k=0}^{\infty}$, obtained by ME-RBK from an initial matrix $X^{0}\in \mathbb{R}^{p\times q}$, in which $X_{:,j}^{0}\in R(A^{T})$, $j=1,2,\cdots,q$ and $(X_{i,:}^{0})^{T}\in R(B)$, $i=1,2,\cdots,p$, converges to the unique least-norm solution $X_{\ast}$ in expectation. Moreover, the solution error satisfies
\begin{align}
E\left[\|X^{k}-X_{\ast}\|_{F}^{2}\right]
\leq\rho^{k}\|X^{0}-X_{\ast}\|_{F}^{2},
\tag{$2.3$}
\label{eq2.3}
\end{align}
where $\rho=1-\frac{2\alpha-\alpha^{2}\|B\|_{2}^{2}}
{\|A\|_{F}^{2}}\sigma_{\min}^{2}(A)\sigma_{\min}^{2}(B)$.
\end{lemma}


\section{The convergence of ME-BK method}

Assuming that the row index $i_{k}$
has been determined at the $k$-th iteration in a cyclic manner, then the $(k+1)$-th estimate vector $X^{k+1}$ of the Kaczmarz method (1.2) is obtained by projecting the current estimate vector $X^{k}$ orthogenically onto the hyperplane $\{x~|~\langle A_{i_{k},:}^{T},x\rangle=b_{i_{k}}\}$.
Inspired by this, we investigate the following ME-BK method for solving (1.1).

\begin{breakablealgorithm}
\caption{~ME-BK method}
{ \begin{algorithmic}[1]
\noindent {\bf Initialization:} Give any initial iterative matrix $X^{0}$.
\State {\bf for} $k=0,1,2,\cdots$ until convergence {\bf do}
\State ~~~~Compute $i_{k}=(k\ mod\ m)+1$
\State ~~~~Compute $X^{k+1}=X^{k}+\frac{\alpha}{\|A_{i_{k},:}\|_{2}^{2}}A_{i_{k},:}^{T}\left(
C_{i_{k},:}-A_{i_{k},:}X^{k}B\right)B^{T}$
\State {\bf end for}
\end{algorithmic}}
\end{breakablealgorithm}

\begin{remark}
The iteration from $k=lm$ to $k=l(m+1)-1$ is called one cycle for $l=0,1,\cdots$. After one cycle, the ME-BK method cyclically sweeps through the $m$ subsystems $A_{i,:}XB=C_{i,:}$, $i=1,2,\cdots,m$ of matrix equation (1.1) and estimation matrix $X^{k}$ is iteratively updated m times.
\end{remark}


Let $X_{\ast}^{0}=X_{\ast}+X^{0}-A^{+}AX^{0}BB^{+}$. 
The convergence of the iterative method is one of the primary considerations to take into account.
We note that Theorem A.3 in [17] served as the inspiration for the following convergence analysis.

\begin{lemma}
The iteration matrix $X^{k}$ obtained by the ME-BK method from an arbitrary initial guess $X^{0}\in \mathbb{R}^{p\times q}$ satisfies
\begin{align*}
X^{k}-A^{+}AX^{k}BB^{+}=X^{0}-A^{+}AX^{0}BB^{+},\tag{$3.1$}
\end{align*}
for $k=0,1,2,\cdots$.
\end{lemma}
\noindent {\bf Proof.}
From the ME-BK method, we obtain
\begin{align*}
&X^{k}=X^{0}+\sum_{s=0}^{k-1}
\frac{\alpha}{\|A_{i_{s},:}\|_{2}^{2}}A_{i_{s},:}^{T}\left(
C_{i_{s},:}-A_{i_{s},:}X^{k}B\right)B^{T}. \tag{$3.2$}
\end{align*}
Subsequently, it is obtained that
\begin{align*}
A^{+}AX^{k}BB^{+}&=A^{+}A\left[X^{0}+\sum_{s=0}^{k-1}
\frac{\alpha}{\|A_{i_{s},:}\|_{2}^{2}}A_{i_{s},:}^{T}\left(
C_{i_{s},:}-A_{i_{s},:}X^{k}B\right)B^{T}\right]BB^{+}\\
&=A^{+}A X^{0}BB^{+}+A^{+}A\left[\sum_{s=0}^{k-1}
\frac{\alpha}{\|A_{i_{s},:}\|_{2}^{2}}A_{i_{s},:}^{T}\left(
C_{i_{s},:}-A_{i_{s},:}X^{k}B\right)B^{T}\right]BB^{+}\\
&=A^{+}A X^{0}BB^{+}+\sum_{s=0}^{k-1}
\frac{\alpha}{\|A_{i_{s},:}\|_{2}^{2}}A_{i_{s},:}^{T}\left(
C_{i_{s},:}-A_{i_{s},:}X^{k}B\right)B^{T}
.\tag{$3.3$}
\end{align*}
The last equality is valid because
\begin{align*}
A^{+}AA_{i_{s},:}^{T}
=(A^{+}A)^{T}(I_{i_{s},:}A)^{T}
=A^{T}(A^{T})^{+}A^{T}I_{i_{s},:}^{T}
=A^{T}I_{i_{s},:}^{T}
=A_{i_{s},:}^{T}
\end{align*}
and
\begin{align*}
B^{T}BB^{+}=B^{T}(BB^{+})^{T}=B^{T}(B^{T})^{+}B^{T}=B^{T}.
\end{align*}
By (3.2) and (3.3), the conclusion is proved. ~~~~~~\fbox {}

\begin{theorem}
Let the matrix equation (1.1) with $A\in \mathbb{R}^{m\times p}$, $B\in \mathbb{R}^{n\times q}$ and $C\in \mathbb{R}^{m\times n}$
be consistent.
Assume that $0<\alpha<\frac{2}{\|B\|_{2}^{2}}$,
then the sequence $\{X^{k}\}_{k=0}^{\infty}$, obtained by the ME-BK method from an arbitrary initial matrix $X^{0}\in \mathbb{R}^{p\times q}$, converges to the solution $X_{\ast}^{0}$.
\end{theorem}
\noindent {\bf Proof.}
By the ME-BK method, one has
\begin{align*}
\|X^{k+1}-X_{\ast}\|_{F}^{2}&=\|X^{k+1}-X^{k}+X^{k}-X_{\ast}\|_{F}^{2}\\
&=\|X^{k+1}-X^{k}\|_{F}^{2}
+2\left\langle X^{k+1}-X^{k}, X^{k}-X_{\ast}\right\rangle
+\|X^{k}-X_{\ast}\|_{F}^{2}\\
&=\frac{\alpha^{2}}{\|A_{i_{k},:}\|_{2}^{4}}
\left\|A_{i_{k},:}^{T}(
C_{i_{k},:}-A_{i_{k},:}X^{k}B)B^{T}\right\|_{F}^{2}\\
&~~+2\frac{\alpha}{\|A_{i_{k},:}\|_{2}^{2}}
\left\langle A_{i_{k},:}^{T}(
C_{i_{k},:}-A_{i_{k},:}X^{k}B)B^{T}, X^{k}-X_{\ast}\right\rangle
+\|X^{k}-X_{\ast}\|_{F}^{2}
.\tag{$3.4$}
\end{align*}
Note that
\begin{align*}
&\frac{\alpha^{2}}{\|A_{i_{k},:}\|_{2}^{4}}
\left\|A_{i_{k},:}^{T}
(C_{i_{k},:}-A_{i_{k},:}X^{k}B)B^{T}\right\|_{F}^{2}\\
&=\frac{\alpha^{2}}{\|A_{i_{k},:}\|_{2}^{4}}
tr[B(C_{i_{k},:}-A_{i_{k},:}X^{k}B)^{T}A_{i_{k},:}
A_{i_{k},:}^{T}(C_{i_{k},:}-A_{i_{k},:}X^{k}B)B^{T}]\\
&=\frac{\alpha^{2}}{\|A_{i_{k},:}\|_{2}^{2}}
tr[B(C_{i_{k},:}-A_{i_{k},:}X^{k}B)^{T}
(C_{i_{k},:}-A_{i_{k},:}X^{k}B)B^{T}]\\
&=\frac{\alpha^{2}}{\|A_{i_{k},:}\|_{2}^{2}}
\left\|
(C_{i_{k},:}-A_{i_{k},:}X^{k}B)B^{T}\right\|_{2}^{2}\\
&\leq \frac{\alpha^{2}\|B\|_{2}^{2}}{\|A_{i_{k},:}\|_{2}^{2}}
\left\|
C_{i_{k},:}-A_{i_{k},:}X^{k}B\right\|_{2}^{2}
.\tag{$3.5$}
\end{align*}
According to $AX_{\ast}B=C$ which can be found in [32], one has
\begin{align*}
&\left\langle A_{i_{k},:}^{T}(C_{i_{k},:}-A_{i_{k},:}X^{k}B)B^{T},
X^{k}-X_{\ast}\right\rangle\\
&=tr[(C_{i_{k},:}-A_{i_{k},:}X^{k}B)^{T}A_{i_{k},:}(X^{k}-X_{\ast})B]\\
&=tr[(C_{i_{k},:}-A_{i_{k},:}X^{k}B)^{T}I_{i_{k},:}A(X^{k}-X_{\ast})B]\\
&=tr[(C_{i_{k},:}-A_{i_{k},:}X^{k}B)^{T}I_{i_{k},:}(AX^{k}B-AX_{\ast}B)]\\
&=tr[(C_{i_{k},:}-A_{i_{k},:}X^{k}B)^{T}(A_{i_{k},:}X^{k}B -C_{i_{k},:})]\\
&=-\left\|
C_{i_{k},:}-A_{i_{k},:}X^{k}B\right\|_{2}^{2}
.\tag{$3.6$}
\end{align*}
Substituting (3.5) and (3.6) into (3.4) yields
\begin{align*}
&\|X^{k+1}-X_{\ast}\|_{F}^{2}-\|X^{k}-X_{\ast}\|_{F}^{2}\\
&\leq \frac{\alpha^{2}\|B\|_{2}^{2}}{\|A_{i_{k},:}\|_{2}^{2}}
\left\|C_{i_{k},:}-A_{i_{k},:}X^{k}B\right\|_{2}^{2}
-\frac{2\alpha}{\|A_{i_{k},:}\|_{2}^{2}}
\left\|C_{i_{k},:}-A_{i_{k},:}X^{k}B\right\|_{2}^{2}\\
&=-\frac{2\alpha-\alpha^{2}\|B\|_{2}^{2}}{\|A_{i_{k},:}\|_{2}^{2}}
\left\|C_{i_{k},:}-A_{i_{k},:}X^{k}B\right\|_{2}^{2}
.\tag{$3.7$}
\end{align*}
If $0<\alpha<\frac{2}{\|B\|_{2}^{2}}$, then $\left\{\|X^{k}-X_{\ast}\|_{F}^{2}\right\}$ is a monotonically decreasing sequence and converges to a limit.
Suppose that the sequence $\left\{\|X^{k}-X_{\ast}\|_{F}^{2}\right\}$ converges to a real number $g\geq 0$.

From (3.7), it has
\begin{align*}
\|X^{(k+1)m}-X_{\ast}\|_{F}^{2}
\leq\|X^{km}-X_{\ast}\|_{F}^{2}
-\sum_{t=1}^{m}\frac{2\alpha-\alpha^{2}\|B\|_{2}^{2}}{\|A_{t,:}\|_{2}^{2}}
\left\|C_{t,:}-A_{t,:}X^{km+t-1}B\right\|_{2}^{2}
.\tag{$3.8$}
\end{align*}
Thus, we have
\begin{align*}
0\leq\sum_{t=1}^{m}\frac{2\alpha-\alpha^{2}\|B\|_{2}^{2}}{\|A_{t,:}\|_{2}^{2}}
\left\|C_{t,:}-A_{t,:}X^{km+t-1}B\right\|_{2}^{2}
\leq \|X^{km}-X_{\ast}\|_{F}^{2}-
\|X^{(k+1)m}-X_{\ast}\|_{F}^{2}
.\tag{$3.9$}
\end{align*}
Obviously, $\lim\limits_{k\rightarrow\infty}\left(\|X^{km}-X_{\ast}\|_{F}^{2}-
\|X^{(k+1)m}-X_{\ast}\|_{F}^{2}\right)
=g-g=0$.
Therefore, it is concluded that
\begin{align*}
\lim\limits_{k\rightarrow\infty}\sum_{t=1}^{m}\frac{2\alpha-\alpha^{2}\|B\|_{2}^{2}}{\|A_{t,:}\|_{2}^{2}}
\left\|C_{t,:}-A_{t,:}X^{km+t-1}B\right\|_{2}^{2}
=0
.\tag{$3.10$}
\end{align*}
This means that
$\lim\limits_{k\rightarrow\infty}
\left\|C_{t,:}-A_{t,:}X^{km+t-1}B\right\|_{2}^{2}
=0$ for $t=1,2,\cdots,m$,
or
\begin{align*}
\lim\limits_{k\rightarrow\infty}
\left(C_{t,:}-A_{t,:}X^{km+t-1}B\right)
=\mathbf{0}, \ t=1,2,\cdots,m
.
\tag{$3.11$}
\end{align*}
It follows from (3.8) that $\{X^{km}-X_{\ast}\}$ is bounded. Hence $\{X^{km}\}$ is bounded.
Suppose subsequence $\{X^{k_{p}m},\ p=1,2,\cdots\}$ converges to a limit $\widehat{X}$. From the ME-BK method and (3.11), we obtain
\begin{align*}
\lim\limits_{p\rightarrow\infty}
\left(X^{k_{p}m+t}-X^{k_{p}m+t-1}\right)
=\mathbf{0}, \ t=1,2,\cdots,m
.
\tag{$3.12$}
\end{align*}
This implies that $\widehat{X}$ is the limit of every subsequences $\{X^{k_{p}m+t-1}\}$ for $t=1,2,\cdots, m$.
From (3.11), we get
\begin{align*}
C_{t,:}=\lim\limits_{p\rightarrow\infty}
A_{t,:}X^{k_{p}m+t-1}B
=A_{t,:}\widehat{X}B
,\ t=1,2,\cdots,m,
\tag{$3.13$}
\end{align*}
which implies that $\widehat{X}$ is a solution of (1.1).

According to (3.7), it is obvious that
\begin{align*}
\|X^{k+1}-\widehat{X}\|_{F}^{2}\leq\|X^{k}-\widehat{X}\|_{F}^{2}.
\tag{$3.14$}
\end{align*}
Therefore,
$\left\{\|X^{k}-\widehat{X}\|_{F}^{2}\right\}$ is a monotonically decreasing and bounded sequence.
Hence, we conclude that $\lim\limits_{k\rightarrow\infty}\|X^{k}-\widehat{X}\|_{F}^{2}$ exists. Because $\lim\limits_{p\rightarrow\infty}\|X^{k_{p}m}-\widehat{X}\|_{F}^{2}=0$, it follows that $\lim\limits_{k\rightarrow\infty}\|X^{k}-\widehat{X}\|_{F}^{2}=0$, i.e., $\lim\limits_{k\rightarrow\infty}X^{k}=\widehat{X}$.

Define $\mathcal{N}:=\{X\in \mathbb{R}^{p\times q}\ |\ AXB=\mathbf{0}\}$.
Let $Z=\widehat{X}-A^{+}A\widehat{X}BB^{+}$.
Because $Z\in\mathcal{N}$ and
\begin{align*}
\langle \widehat{X}-Z,Z\rangle
&=\langle A^{+}A\widehat{X}BB^{+},\widehat{X}-A^{+}A\widehat{X}BB^{+}\rangle\\
&=\langle A^{+}A\widehat{X}BB^{+},\widehat{X}\rangle
-\langle A^{+}A\widehat{X}BB^{+},A^{+}A\widehat{X}BB^{+}\rangle\\
&=tr(\widehat{X}^{T}A^{+}A\widehat{X}BB^{+})
-tr((BB^{+})^{T}\widehat{X}^{T}(A^{+}A)^{T}A^{+}A\widehat{X}BB^{+})\\
&=tr(\widehat{X}^{T}A^{+}A\widehat{X}BB^{+})
-tr(\widehat{X}^{T}A^{+}AA^{+}A\widehat{X}BB^{+}BB^{+})\\
&=tr(\widehat{X}^{T}A^{+}A\widehat{X}BB^{+})
-tr(\widehat{X}^{T}A^{+}A\widehat{X}BB^{+})\\
&=0,\tag{$3.15$}
\end{align*}
it has $\|\widehat{X}\|_{F}^{2}=\|\widehat{X}-Z\|_{F}^{2}+\|Z\|_{F}^{2}\geq\|\widehat{X}-Z\|_{F}^{2}$.
Therefore, $\widehat{X}-Z$ is the minimum norm solution, i.e. $\widehat{X}-Z=X_{\ast}$.
According to Lemma 3.2, we get
\begin{align*}
\widehat{X}-A^{+}A\widehat{X}BB^{+}
&=\lim\limits_{k\rightarrow\infty}(X^{k}-A^{+}AX^{k}BB^{+})\\
&=\lim\limits_{k\rightarrow\infty}(X^{0}-A^{+}AX^{0}BB^{+})\\
&=X^{0}-A^{+}AX^{0}BB^{+}.\tag{$3.16$}
\end{align*}
It follows from (3.16) that $\widehat{X}=X_{\ast}+X^{0}-A^{+}AX^{0}BB^{+}$. The proof is complete.~~~~~~\fbox {}

\section{ME-GRBK method and its relaxation and deterministic versions}



According to the method outlined in the literature [19], if the norm of the $i$-th row of the residual matrix is greater than that of the $j$-th row at the $k$-th iteration,
we select the $i$-th row with greater probability, thereby eliminating the larger residual matrix entries as preferentially as possible.
In this way, we derive the following ME-GRBK method for solving (1.1).

\begin{breakablealgorithm}
\caption{~ME-GRBK method}
{ \begin{algorithmic}[1]
\noindent{\bf Initialization:} Give any initial iterative matrix $X^{0}$ and set $R^{0}=C-AX^{0}B$.
\State {\bf for} $k=0,1,2,\cdots$ until convergence {\bf do}
\State ~~~~Compute $\zeta_{k}=\frac{1}{2}\left[
\frac{1}{\|R^{k}\|_{F}^{2}}
\max\limits_{1\leq i\leq m}\left(\frac{\|R_{i,:}^{k}\|_{2}^{2}}
{\|A_{i,:}\|_{2}^{2}}\right)
+\frac{1}{\|A\|_{F}^{2}}\right]$
\State ~~~~Determine the index set
\begin{align*}
\mathcal{J}_{k}=\left\{i_{k}~|~
\|R_{i_{k},:}^{k}\|_{2}^{2}\geq\zeta_{k}
\|A_{i_{k},:}\|_{2}^{2}\|R^{k}\|_{F}^{2}
\right\}
\end{align*}
\State ~~~~Compute the $i$th row $\widetilde{R}_{i,:}^{k}$ of the matrix $\widetilde{R}^{k}$ according to
\begin{align*}
\widetilde{R}_{i,:}^{k}=
\left\{
             \begin{array}{ll}
               R_{i,:}^{k}, & if~ i\in\mathcal{J}_{k}\\
               0, &otherwise\\
             \end {array}
           \right.
\end{align*}
\State ~~~~Choose $i_{k}\in\mathcal{J}_{k}$ with probability
$ {\rm Pr}({\rm row=}i_{k})
=\frac{\|\widetilde{R}_{i_{k},:}^{k}\|_{2}^{2}}{\|\widetilde{R}^{k}\|_{F}^{2}}$
\State ~~~~Update $X^{k+1}=X^{k}+\frac{\alpha}{\|A_{i_{k},:}\|_{2}^{2}}
A_{i_{k},:}^{T}
R_{i_{k},:}^{k}B^{T}$
\State ~~~~Compute  $R^{k+1}=R^{k}-\frac{\alpha}{\|A_{i_{k},:}\|_{2}^{2}}
(AA^{T})_{:,i_{k}}R_{i_{k},:}^{k}B^{T}B$
\State {\bf end for}
\end{algorithmic}}
\end{breakablealgorithm}

In practice, we first compute $AA^{T}$ and $B^{T}B$ and bring them directly into the update of $R^{k+1}$ during the iteration of the ME-GRBK method. This eliminates the requirement to recalculate $AA^{T}$ and $B^{T}B$ at each iteration, saving computation and speeding up the convergence of the algorithm.

\begin{remark}
The index set $\mathcal{J}_{k}$ is updated once per iteration and is nonempty for all iteration numbers $k$. The reason is as follows.
It is obvious that
\begin{align*}
\max\limits_{1\leq i\leq m}\left(\frac{\|R_{i,:}^{k}\|_{2}^{2}}
{\|A_{i,:}\|_{2}^{2}}\right)
\geq\sum_{i=1}^{m}
\frac{\|A_{i,:}\|_{2}^{2}}{\|A\|_{F}^{2}}
\frac{\|R_{i,:}^{k}\|_{2}^{2}}
{\|A_{i,:}\|_{2}^{2}}
=\frac{\|R^{k}\|_{F}^{2}}{\|A\|_{F}^{2}}.
\end{align*}
Let $\frac{\|R_{s,:}^{k}\|_{2}^{2}}
{\|A_{s,:}\|_{2}^{2}}=\max\limits_{1\leq i\leq m}\left(\frac{\|R_{i,:}^{k}\|_{2}^{2}}
{\|A_{i,:}\|_{2}^{2}}\right)$, then we obtain
\begin{align*}
\frac{\|R_{s,:}^{k}\|_{2}^{2}}
{\|A_{s,:}\|_{2}^{2}}=\max\limits_{1\leq i\leq m}\left(\frac{\|R_{i,:}^{k}\|_{2}^{2}}
{\|A_{i,:}\|_{2}^{2}}\right)
\geq\frac{1}{2}\max\limits_{1\leq i\leq m}\left(\frac{\|R_{i,:}^{k}\|_{2}^{2}}
{\|A_{i,:}\|_{2}^{2}}\right)
+\frac{1}{2}\frac{\|R^{k}\|_{F}^{2}}{\|A\|_{F}^{2}}
=\eta_{k}
\|R^{k}\|_{F}^{2}.
\end{align*}
Hence, $s\in\mathcal{J}_{k}$ and $\mathcal{J}_{k}\neq\emptyset$.
\end{remark}

We provide the convergence of the ME-GRBK method in the following theorem, which is inspired by [11, 19].

\begin{theorem}
Let the matrix equation (1.1) with $A\in \mathbb{R}^{m\times p}$, $B\in \mathbb{R}^{n\times q}$ and $C\in \mathbb{R}^{m\times n}$
be consistent.
Assume that $0<\alpha<\frac{2}{\|B\|_{2}^{2}}$,
then the sequence $\{X^{k}\}_{k=0}^{\infty}$, obtained by the ME-GRBK method from an initial matrix $X^{0}\in \mathbb{R}^{p\times q}$ with $vec(X^{0})\in R(B\otimes A^{T})$, converges to the unique least-norm solution $X_{\ast}$ in expectation.
Moreover, the solution error
obeys
\begin{align*}
E\left[\|X^{1}-X_{\ast}\|_{F}^{2}\right]
\leq\left(1-\frac{(2\alpha-\alpha^{2}\|B\|_{2}^{2})\|B\|_{F}^{2}}{\kappa_{F,2}^{2}(B^{T}\otimes A)}\right)
\|X^{0}-X_{\ast}\|_{F}^{2},\tag{$4.1$}
\end{align*}
and
\begin{align*}
E_{k}\left[\|X^{k+1}-X_{\ast}\|_{F}^{2}\right]
\leq\left(1-\frac{(2\alpha-\alpha^{2}\|B\|_{2}^{2})\varphi_{k}}{\kappa_{F,2}^{2}(B^{T}\otimes A)}\right)
\|X^{k}-X_{\ast}\|_{F}^{2},\tag{$4.2$}
\end{align*}
where $\varphi_{k}=\frac{\|B^{T}\otimes A\|_{F}^{2}}{2\gamma_{k}}
+\frac{\|B\|_{F}^{2}}{2}$, $\gamma_{k}=\|A\|_{F}^{2}-\sum_{i\in\Omega_{k}}
\|A_{i,:}\|_{2}^{2}$ and $\Omega_{k}=\left\{i~|~\left\|R^{k}_{i,:}\right\|_{2}^{2}=0, i\in[m]\right\}$.
As a result, it has
\begin{align*}
E\left[\|X^{k+1}-X_{\ast}\|_{F}^{2}\right]
\leq\displaystyle\prod_{l=0}^{k}\rho_{l}
\|X^{0}-X_{\ast}\|_{F}^{2},\tag{$4.3$}
\end{align*}
where $\rho_{0}=1-\frac{(2\alpha-\alpha^{2}\|B\|_{2}^{2})\|B\|_{F}^{2}}{\kappa_{F,2}^{2}(B^{T}\otimes A)}$ and $\rho_{l}=1-\frac{(2\alpha-\alpha^{2}\|B\|_{2}^{2})\varphi_{l}}{\kappa_{F,2}^{2}(B^{T}\otimes A)}$ for $l=1,2,\cdots,k$.
\end{theorem}
\noindent {\bf Proof.}
It follows from (3.7) that
\begin{align*}
\|X^{k+1}-X_{\ast}\|_{F}^{2}
\leq\|X^{k}- X_{\ast}\|_{F}^{2} -\frac{2\alpha-\alpha^{2}\|B\|_{2}^{2}}{\|A_{i_{k},:}\|_{2}^{2}}
\left\|C_{i_{k},:}-A_{i_{k},:}X^{k}B\right\|_{2}^{2}
.\tag{$4.4$}
\end{align*}
According to the expected value conditional on the first $k$ iterations, it has
\begin{align*}
&E_{k}\left[\frac{2\alpha-\alpha^{2}\|B\|_{2}^{2}}{\|A_{i_{k},:}\|_{2}^{2}}
\left\|C_{i_{k},:}-A_{i_{k},:}X^{k}B\right\|_{2}^{2}\right]\\
&=(2\alpha-\alpha^{2}\|B\|_{2}^{2})
\sum_{i_{k}=1}^{m}\frac{\|\widetilde{R}_{i_{k},:}^{k}\|_{2}^{2}}{\|\widetilde{R}^{k}\|_{F}^{2}}
\frac{\left\|C_{i_{k},:}-A_{i_{k},:}X^{k}B\right\|_{2}^{2}}
{\|A_{i_{k},:}\|_{2}^{2}}\\
&\geq(2\alpha-\alpha^{2}\|B\|_{2}^{2})
\sum_{i_{k}\in\mathcal{J}_{k}}
\frac{\|\widetilde{R}_{i_{k},:}^{k}\|_{2}^{2}}
{\sum\limits_{i_{k}\in\mathcal{J}_{k}}\|\widetilde{R}_{i_{k},:}^{k}\|_{2}^{2}}
\zeta_{k}
\|R^{k}\|_{F}^{2}\\
&=(2\alpha-\alpha^{2}\|B\|_{2}^{2})\zeta_{k}
\|A(X^{k}- X_{\ast})B\|_{F}^{2}\\
&=(2\alpha-\alpha^{2}\|B\|_{2}^{2})\zeta_{k}
\|(B^{T}\otimes A)vec(X^{k}- X_{\ast})\|_{2}^{2}
.\tag{$4.5$}
\end{align*}
We note that $vec(X^{0})\in R(B\otimes A^{T})$ and
\begin{align*}
vec(X_{\ast})
&=((B^{+})^{T}\otimes A^{+})vec(C)\\
&=\left(\left(BB^{+}(B^{+})^{T}\right)\otimes \left(A^{T}(A^{+})^{T}A^{+}\right)\right)vec(C)\\
&=(B\otimes A^{T})\left(\left(B^{+}(B^{+})^{T}\right)\otimes \left((A^{+})^{T}A^{+}\right)\right)vec(C)\in R(B\otimes A^{T})
,
\end{align*}
so it has $vec(X^{0}-X_{\ast})\in R(B\otimes A^{T})$.
We also have
\begin{align*}
vec\left(\frac{\alpha}{\|A_{i_{k-1},:}\|_{2}^{2}}
A_{i_{k-1},:}^{T}
R_{i_{k-1},:}^{k-1}B^{T}\right)
&=\frac{\alpha}{\|A_{i_{k-1},:}\|_{2}^{2}}(B\otimes A_{i_{k-1},:}^{T})
vec(R_{i_{k-1},:}^{k-1})\\
&=\frac{\alpha}{\|A_{i_{k-1},:}\|_{2}^{2}}(B\otimes A^{T})(I\otimes I_{i_{k-1},:}^{T})
vec(R_{i_{k-1},:}^{k-1})\in R(B\otimes A^{T})
.
\end{align*}
Because $X^{k}-X_{\ast}=X^{k-1}-X_{\ast}
+\frac{\alpha}{\|A_{i_{k-1},:}\|_{2}^{2}}
A_{i_{k-1},:}^{T}
R_{i_{k-1},:}^{k-1}B^{T}$, we get $vec(X^{k}-X_{\ast})\in R(B\otimes A^{T})$ by induction. From Lemma 2.1 and $0<\alpha<\frac{2}{\|B\|_{2}^{2}}$, it holds
\begin{align*}
E_{k}\left[\frac{2\alpha-\alpha^{2}\|B\|_{2}^{2}}{\|A_{i_{k},:}\|_{2}^{2}}
\left\|C_{i_{k},:}\!-\!A_{i_{k},:}X^{k}B\right\|_{2}^{2}\right]
&\geq(2\alpha-\alpha^{2}\|B\|_{2}^{2})\zeta_{k}
\frac{\|B^{T}\otimes A\|_{F}^{2}}
{\kappa_{F,2}^{2}(B^{T}\otimes A)}
\|X^{k}- X_{\ast}\|_{F}^{2}
.\tag{$4.6$}
\end{align*}
By (4.4) and (4.6), we obtain
\begin{align*}
E_{k}\left[\|X^{k+1}-X_{\ast}\|_{F}^{2}\right]
\leq\|X^{k}- X_{\ast}\|_{F}^{2} -(2\alpha-\alpha^{2}\|B\|_{2}^{2})\zeta_{k}
\frac{\|B^{T}\otimes A\|_{F}^{2}}
{\kappa_{F,2}^{2}(B^{T}\otimes A)}
\|X^{k}- X_{\ast}\|_{F}^{2}
.\tag{$4.7$}
\end{align*}
In addition, we have
\begin{align*}
\zeta_{0}\|A\|_{F}^{2}
&=\frac{1}{2}\frac{\max\limits_{1\leq i\leq m}\left(\frac{\|R_{i,:}^{0}\|_{2}^{2}}
{\|A_{i,:}\|_{2}^{2}}\right)}
{\frac{\|R^{0}\|_{F}^{2}}
{\|A\|_{F}^{2}}}
+\frac{1}{2}\\
&=\frac{1}{2}\frac{\max\limits_{1\leq i\leq m}\left(\frac{\|R_{i,:}^{0}\|_{2}^{2}}
{\|A_{i,:}\|_{2}^{2}}\right)}
{\sum_{i\in[m]}
\frac{\|A_{i,:}\|_{2}^{2}}
{\|A\|_{F}^{2}}
\frac{\|R_{i,:}^{0}\|_{2}^{2}}
{\|A_{i,:}\|_{2}^{2}}}
+\frac{1}{2}\\
&\geq\frac{1}{2}\frac{1}
{\sum_{i\in[m]}
\frac{\|A_{i,:}\|_{2}^{2}}
{\|A\|_{F}^{2}}}
+\frac{1}{2}\\
&=1,
\tag{$4.8$}
\end{align*}
and
\begin{align*}
\zeta_{k}\|A\|_{F}^{2}
&=\frac{1}{2}\frac{\max\limits_{1\leq i\leq m}\left(\frac{\|R_{i,:}^{k}\|_{2}^{2}}
{\|A_{i,:}\|_{2}^{2}}\right)}
{\frac{\|R^{k}\|_{F}^{2}}
{\|A\|_{F}^{2}}}
+\frac{1}{2}\\
&=\frac{1}{2}\frac{\max\limits_{1\leq i\leq m}\left(\frac{\|R_{i,:}^{k}\|_{2}^{2}}
{\|A_{i,:}\|_{2}^{2}}\right)}
{\sum_{i\in[m]}
\frac{\|A_{i,:}\|_{2}^{2}}
{\|A\|_{F}^{2}}
\frac{\|R_{i,:}^{k}\|_{2}^{2}}
{\|A_{i,:}\|_{2}^{2}}}
+\frac{1}{2}\\
&=\frac{1}{2}\frac{\max\limits_{1\leq i\leq m}\left(\frac{\|R_{i,:}^{k}\|_{2}^{2}}
{\|A_{i,:}\|_{2}^{2}}\right)}
{\left(\sum_{i\in[m]}-\sum_{i\in\Omega_{k}}\right)
\frac{\|A_{i,:}\|_{2}^{2}}
{\|A\|_{F}^{2}}
\frac{\|R_{i,:}^{k}\|_{2}^{2}}
{\|A_{i,:}\|_{2}^{2}}}
+\frac{1}{2}\\
&\geq\frac{1}{2}\frac{1}
{\left(\sum_{i\in[m]}-\sum_{i\in\Omega_{k}}\right)
\frac{\|A_{i,:}\|_{2}^{2}}
{\|A\|_{F}^{2}}}
+\frac{1}{2}
.\tag{$4.9$}
\end{align*}
From (4.8) and (4.9) it follows that
\begin{align*}
\zeta_{0}
\geq\frac{1}{\|A\|_{F}^{2}},
\tag{$4.10$}
\end{align*}
and
\begin{align*}
\zeta_{k}
\geq\frac{1}{2}\left(\frac{1}
{\left(\sum_{i\in[m]}-\sum_{i\in\Omega_{k}}\right)
\|A_{i,:}\|_{2}^{2}}
+\frac{1}{\|A\|_{F}^{2}}\right).
\tag{$4.11$}
\end{align*}
Noting that $\|B^{T}\otimes A\|_{F}=\| A\|_{F}\|B\|_{F}$, then the estimates in (4.1) and (4.2) are derived from (4.7), (4.10) and (4.11).
Taking the full expectation on both sides of inequation (4.2) gives
\begin{align*}
E\left[\|X^{k+1}-X_{\ast}\|_{F}^{2}\right]
\leq\left(1-\frac{(2\alpha-\alpha^{2}\|B\|_{2}^{2})\varphi_{k}}{\kappa_{F,2}^{2}(B^{T}\otimes A)}\right)
E\left[\|X^{k}-X_{\ast}\|_{F}^{2}\right], k=1,2,\cdots.\tag{$4.12$}
\end{align*}
An induction on the iteration number $k$ gives us the estimate (4.3).
~~~~~~\fbox {}

\begin{remark}
Since $\frac{1}
{2\left(\|A\|_{F}^{2}-\sum_{i\in\Omega_{k}}
\|A_{i,:}\|_{2}^{2}\right)}
+\frac{1}{2\|A\|_{F}^{2}}\geq\frac{1}{\|A\|_{F}^{2}}$, this gives
\begin{align*}
\rho_{k}=1-\frac{(2\alpha-\alpha^{2}\|B\|_{2}^{2})\varphi_{k}}{\kappa_{F,2}^{2}(B^{T}\otimes A)}\leq1-\frac{2\alpha-\alpha^{2}\|B\|_{2}^{2}}
{\|A\|_{F}^{2}}\sigma_{\min}^{2}(A)\sigma_{\min}^{2}(B)=\rho.
\end{align*}
As a result, the convergence factor of the ME-GRBK method is smaller than that of the ME-RBK method in literature [12].
\end{remark}

Motivated by the work in literature [20], we further generalise the ME-GRBK method by introducing a relaxation factor $\theta\in(0,1)$ in the probability criterion, and obtain a relaxed version of the ME-GRBK method for solving (1.1), which is presented as follows.

\begin{breakablealgorithm}
\caption{~ME-RGRBK method}
{ \begin{algorithmic}[1]
\noindent{\bf Initialization:} Give any initial iterative matrix $X^{0}$ and set $R^{0}=C-AX^{0}B$. Give a relaxation factor $\theta\in(0,1)$.
\State {\bf for} $k=0,1,2,\cdots$ until convergence {\bf do}
\State ~~~~Compute $\xi_{k}=
\frac{\theta}{\|R^{k}\|_{F}^{2}}
\max\limits_{1\leq i\leq m}\left(\frac{\|R_{i,:}^{k}\|_{2}^{2}}
{\|A_{i,:}\|_{2}^{2}}\right)
+\frac{1-\theta}{\|A\|_{F}^{2}}$
\State ~~~~Determine the index set
\begin{align*}
\mathcal{H}_{k}=\left\{i_{k}~|~
\|R_{i_{k},:}^{k}\|_{2}^{2}\geq\xi_{k}
\|A_{i_{k},:}\|_{2}^{2}\|R^{k}\|_{F}^{2}
\right\}
\end{align*}
\State ~~~~Compute the $i$th row $\widetilde{R}_{i,:}^{k}$ of the matrix $\widetilde{R}^{k}$ according to
\begin{align*}
\widetilde{R}_{i,:}^{k}=
\left\{
             \begin{array}{ll}
               R_{i,:}^{k}, & if~ i\in\mathcal{H}_{k}\\
               0, &otherwise\\
             \end {array}
           \right.
\end{align*}
\State ~~~~Choose $i_{k}\in\mathcal{H}_{k}$ with probability
$ {\rm Pr}({\rm row=}i_{k})
=\frac{\|\widetilde{R}_{i_{k},:}^{k}\|_{2}^{2}}{\|\widetilde{R}^{k}\|_{F}^{2}}$
\State ~~~~Update $X^{k+1}=X^{k}+\frac{\alpha}{\|A_{i_{k},:}\|_{2}^{2}}
A_{i_{k},:}^{T}
R_{i_{k},:}^{k}B^{T}$
\State ~~~~Compute  $R^{k+1}=R^{k}-\frac{\alpha}{\|A_{i_{k},:}\|_{2}^{2}}
(AA^{T})_{:,i_{k}}R_{i_{k},:}^{k}B^{T}B$
\State {\bf end for}
\end{algorithmic}}
\end{breakablealgorithm}

\begin{remark}
The most significant improvement of the ME-RGRBK method over the ME-GRBK method is the introduction of a relaxation factor $\theta$ in the index set, which allows for a more flexible and adaptive selection of the appropriate index set. In particular, the ME-RGRBK method is converted to the ME-GRBK method when $\theta$ takes $\frac{1}{2}$.
And the reason for the non-emptiness of the set $\mathcal{H}_{k}$ is analogous to the content of Remark 4.1.
\end{remark}

\begin{theorem}
Let the matrix equation (1.1) with $A\in \mathbb{R}^{m\times p}$, $B\in \mathbb{R}^{n\times q}$ and $C\in \mathbb{R}^{m\times n}$
be consistent.
Assume that $0<\alpha<\frac{2}{\|B\|_{2}^{2}}$,
then the sequence $\{X^{k}\}_{k=0}^{\infty}$, obtained by the ME-RGRBK method from an initial matrix $X^{0}\in \mathbb{R}^{p\times q}$ with $vec(X^{0})\in R(B\otimes A^{T})$, converges to the unique least-norm solution $X_{\ast}$ in expectation.
Moreover, the solution error in expectation for the iteration sequence $\{X^{k}\}_{k=0}^{\infty}$ satisfies
\begin{align*}
E\left[\|X^{1}-X_{\ast}\|_{F}^{2}\right]
\leq\left(1-\frac{(2\alpha-\alpha^{2}\|B\|_{2}^{2})\|B\|_{F}^{2}}{\kappa_{F,2}^{2}(B^{T}\otimes A)}\right)
\|X^{0}-X_{\ast}\|_{F}^{2},\tag{$4.13$}
\end{align*}
and
\begin{align*}
E_{k}\left[\|X^{k+1}-X_{\ast}\|_{F}^{2}\right]
\leq\left(1-\frac{(2\alpha-\alpha^{2}\|B\|_{2}^{2})\phi_{k}}{\kappa_{F,2}^{2}(B^{T}\otimes A)}\right)
\|X^{k}-X_{\ast}\|_{F}^{2},\tag{$4.14$}
\end{align*}
where $\phi_{k}=\frac{\theta\|B^{T}\otimes A\|_{F}^{2}}{\gamma_{k}}
+(1-\theta)\|B\|_{F}^{2}$, $\gamma_{k}=\|A\|_{F}^{2}-\sum_{i\in\Omega_{k}}
\|A_{i,:}\|_{2}^{2}$ and $\Omega_{k}=\left\{i\,|\,\left\|R^{k}_{i,:}\right\|_{2}^{2}=0, i\in[m]\right\}$.
As a result, it has
\begin{align*}
E\left[\|X^{k+1}-X_{\ast}\|_{F}^{2}\right]
\leq\displaystyle\prod_{l=0}^{k}\rho_{l,\theta}
\|X^{0}-X_{\ast}\|_{F}^{2},\tag{$4.15$}
\end{align*}
where $\rho_{0,\theta}=\rho_{0}=1-\frac{(2\alpha-\alpha^{2}\|B\|_{2}^{2})\|B\|_{F}^{2}}{\kappa_{F,2}^{2}(B^{T}\otimes A)}$ and $\rho_{l,\theta}=1-\frac{(2\alpha-\alpha^{2}\|B\|_{2}^{2})\phi_{l}}{\kappa_{F,2}^{2}(B^{T}\otimes A)}$ for $l=1,2,\cdots,k$.
\end{theorem}
\noindent {\bf Proof.} The proof is similar to Theorem 4.2 and is omitted here.
~~~~~~\fbox {}\

\begin{remark}
It can be obtained from $\frac{\theta}
{\left(\|A\|_{F}^{2}-\sum_{i\in\Omega_{k}}
\|A_{i,:}\|_{2}^{2}\right)}
+\frac{1-\theta}{\|A\|_{F}^{2}}\geq\frac{1}{\|A\|_{F}^{2}}$ that
\begin{align*}
\rho_{k,\theta}=1-\frac{(2\alpha-\alpha^{2}\|B\|_{2}^{2})\phi_{k}}{\kappa_{F,2}^{2}(B^{T}\otimes A)}\leq1-\frac{2\alpha-\alpha^{2}\|B\|_{2}^{2}}
{\|A\|_{F}^{2}}\sigma_{\min}^{2}(A)\sigma_{\min}^{2}(B)=\rho.
\end{align*}
Thus, the convergence factor of the ME-RGRBK method is smaller than that of the ME-RBK method [12].
\end{remark}

We note that $\frac{\partial\rho_{k,\theta}}{\partial\theta}
=\frac{(2\alpha-\alpha^{2}\|B\|_{2}^{2})
\|B^{T}\otimes A\|_{F}^{2}}{\kappa_{F,2}^{2}(B^{T}\otimes A)}
\left(\frac{1}{\|A\|_{F}^{2}}-\frac{1}{\|A\|_{F}^{2}-\sum_{i\in\Omega_{k}}
\|A_{i,:}\|_{2}^{2}}\right)<0$, so $\rho_{k,\theta}$ is monotonically decreasing with respect to $\theta$ for $\theta\in(0,1)$. Therefore, it has $\rho_{k,\theta}>\rho_{k}$ if $\theta\in(0,\frac{1}{2})$ and $\rho_{k,\theta}<\rho_{k}$ if $\theta\in(\frac{1}{2},1)$. In addition, the bound on the convergence factor of ME-RGRBK method is minimised when $\theta=1$. On the basis of this finding, the following MR-MWRBK method is derived to solve (1. 1), which is a deterministic version of the ME-GRBK method.

\begin{breakablealgorithm}
\caption{~ME-MWRBK method}
{ \begin{algorithmic}[1]
\noindent{\bf Initialization:} Give any initial iterative matrix $X^{0}$ and set $R^{0}=C-AX^{0}B$.
\State {\bf for} $k=0,1,2,\cdots$ until convergence {\bf do}
\State ~~~~Determine the index set
\begin{align*}
\mathcal{G}_{k}=\left\{i_{k}~\bigg|~
\frac{\|R_{i_{k},:}^{k}\|_{2}^{2}}
{\|A_{i_{k},:}\|_{2}^{2}}
=\max\limits_{1\leq i\leq m}\left(\frac{\|R_{i,:}^{k}\|_{2}^{2}}
{\|A_{i,:}\|_{2}^{2}}\right)
\right\}
\end{align*}
\State ~~~~Compute $X^{k+1}=X^{k}+\frac{\alpha}{\|A_{i_{k},:}\|_{2}^{2}}
A_{i_{k},:}^{T}\left(
C_{i_{k},:}-A_{i_{k},:}X^{k}B\right)B^{T}$, $i_{k}\in\mathcal{G}_{k}$
\State {\bf end for}
\end{algorithmic}}
\end{breakablealgorithm}

\begin{remark}
In the event that multiple alternatives are presented in $\mathcal{G}_{k}$, it is sufficient to select any one of them.
\end{remark}

A theoretical estimate of the convergence rate of the ME-MWRBK method is given below, which is along the same lines as in [11, 25].

\begin{theorem}
Let the matrix equation (1.1) with $A\in \mathbb{R}^{m\times p}$, $B\in \mathbb{R}^{n\times q}$ and $C\in \mathbb{R}^{m\times n}$
be consistent.
Assume that $0<\alpha<\frac{2}{\|B\|_{2}^{2}}$,
then the sequence $\{X^{k}\}_{k=0}^{\infty}$, obtained by the ME-MWRBK method from an initial matrix $X^{0}\in \mathbb{R}^{p\times q}$ with $vec(X^{0})\in R(B\otimes A^{T})$, converges to the unique least-norm solution $X_{\ast}$.
Moreover, the squared error satisfies
\begin{align*}
\|X^{k}-X_{\ast}\|_{F}^{2}
\leq\displaystyle\prod_{l=0}^{k}\widetilde{\rho_{l}}
\|X^{0}-X_{\ast}\|_{F}^{2},\tag{$4.16$}
\end{align*}
where $\widetilde{\rho}_{0}=1-\frac{(2\alpha-\alpha^{2}\|B\|_{2}^{2})\|B\|_{F}^{2}}{\kappa_{F,2}^{2}(B^{T}\otimes A)}$ and $\rho_{l}=1-\frac{(2\alpha-\alpha^{2}\|B\|_{2}^{2})}
{\gamma_{k}}\sigma_{\min}^{2}(A)\sigma_{\min}^{2}(B)$ for $l=1,2,\cdots,k$ with $\gamma_{k}=\|A\|_{F}^{2}-\sum_{i\in\Omega_{k}}
\|A_{i,:}\|_{2}^{2}$ and $\Omega_{k}=\left\{i\,|\,\left\|R^{k}_{i,:}\right\|_{2}^{2}=0, i\in[m]\right\}$.
\end{theorem}
\noindent {\bf Proof.}
For $k=0$, it holds
\begin{align*}
\|X^{1}-X_{\ast}\|_{F}^{2}
&\leq\|X^{0}- X_{\ast}\|_{F}^{2} -\frac{2\alpha-\alpha^{2}\|B\|_{2}^{2}}{\|A_{i_{0},:}\|_{2}^{2}}
\left\|C_{i_{0},:}-A_{i_{0},:}X^{0}B\right\|_{2}^{2}\\
&=\|X^{0}- X_{\ast}\|_{F}^{2} -(2\alpha-\alpha^{2}\|B\|_{2}^{2})
\frac{\left\|C_{i_{0},:}-A_{i_{0},:}X^{0}B\right\|_{2}^{2}}
{\|A_{i_{0},:}\|_{2}^{2}}
\frac{\left\|C-AX^{0}B\right\|_{F}^{2}}
{\sum\limits_{i=1}^{m}\frac{\left\|C_{i,:}-A_{i,:}X^{0}B\right\|_{2}^{2}}
{\|A_{i,:}\|_{2}^{2}}\|A_{i,:}\|_{2}^{2}}\\
&\leq\|X^{0}- X_{\ast}\|_{F}^{2} -(2\alpha-\alpha^{2}\|B\|_{2}^{2})
\frac{\left\|A(X^{0}- X_{\ast})B\right\|_{F}^{2}}
{\sum\limits_{i=1}^{m}\|A_{i,:}\|_{2}^{2}}\\
&\leq\left(1-\frac{(2\alpha-\alpha^{2}\|B\|_{2}^{2})\|B\|_{F}^{2}}{\kappa_{F,2}^{2}(B^{T}\otimes A)}\right)
\|X^{0}-X_{\ast}\|_{F}^{2}
.\tag{$4.17$}
\end{align*}
For $k\geq1$, we have
\begin{align*}
\|X^{k+1}-X_{\ast}\|_{F}^{2}
&\leq\|X^{k}- X_{\ast}\|_{F}^{2} -\frac{2\alpha-\alpha^{2}\|B\|_{2}^{2}}{\|A_{i_{k},:}\|_{2}^{2}}
\left\|C_{i_{k},:}-A_{i_{k},:}X^{k}B\right\|_{2}^{2}\\
&=\|X^{k}- X_{\ast}\|_{F}^{2} -(2\alpha-\alpha^{2}\|B\|_{2}^{2})
\frac{\left\|R_{i_{k},:}^{k}\right\|_{2}^{2}}
{\|A_{i_{k},:}\|_{2}^{2}}
\frac{\left\|R^{k}\right\|_{F}^{2}}
{\sum\limits_{i=1}^{m}\frac{\left\|R_{i,:}^{k}\right\|_{2}^{2}}
{\|A_{i,:}\|_{2}^{2}}\|A_{i,:}\|_{2}^{2}}\\
&=\|X^{k}- X_{\ast}\|_{F}^{2} -(2\alpha-\alpha^{2}\|B\|_{2}^{2})
\frac{\left\|R_{i_{k},:}^{k}\right\|_{2}^{2}}
{\|A_{i_{k},:}\|_{2}^{2}}
\frac{\left\|R^{k}\right\|_{F}^{2}}
{\left(\sum_{i\in[m]}-\sum_{i\in\Omega_{k}}\right)
\frac{\left\|R_{i,:}^{k}\right\|_{2}^{2}}
{\|A_{i,:}\|_{2}^{2}}\|A_{i,:}\|_{2}^{2}}\\
&\leq\|X^{k}- X_{\ast}\|_{F}^{2} -(2\alpha-\alpha^{2}\|B\|_{2}^{2})
\frac{\left\|A(X^{k}- X_{\ast})B\right\|_{F}^{2}}
{\left(\sum_{i\in[m]}-\sum_{i\in\Omega_{k}}\right)
\|A_{i,:}\|_{2}^{2}}\\
&\leq\left(1-\frac{(2\alpha-\alpha^{2}\|B\|_{2}^{2})}
{\gamma_{k}}\sigma_{\min}^{2}(A)\sigma_{\min}^{2}(B)\right)
\|X^{k}-X_{\ast}\|_{F}^{2}
.\tag{$4.18$}
\end{align*}
The reason why the last line of (4.17) and (4.18) holds is due to $vec(X^{k}-X_{\ast})\in R(B\otimes A^{T})$ and Lemma 2.1. Induction on the iteration number $k$ gives (4.16).
~~~~~~\fbox {}

\section{Numerical examples}

This section provides some numerical examples to validate the primary theoretical findings.
All experiments are tested via MATLAB R2021a on a desktop PC with
Intel(R) Core(TM) i7-9700K CPU @ 3.60 GHz.

We examine the methods presented in this paper from the aspects of the number of iteration step (denoted by `IT'), elapsed CPU time in seconds (denoted by `CPU') and the relative residual norm (denoted by `RRN').

For each example, the matrix $X$ is produced by $\rm randn(p,q)$ and $\alpha$ is taken to be $\frac{1}{\|B\|_{2}^{2}}$.
The following matrices $A$ and $B$ include sparse matrices from the University of Florida sparse matrix collection [33] and dense (sparse) matrices randomly generated by the function randn (sprandn).
To create a consistent system, the matrix $C$ is generated by $C=AXB$.
All computations are terminated once
$\rm{RRN}\leq 10^{-6}$.

For the tested matrices, we show their size, rank and density, where the density is denoted as
\begin{align*}
\text{density}:=
\frac{\text{number of nonzeros of an}\ m\times n\ \text{matrix}}
{m\times n}.
\end{align*}

\noindent {\bf Example 5.1.} In this example, we mainly verify the conclusion in Theorem 3.3. We choose six sets of matrices for testing the convergence of the ME-BK method, as shown in Table 1. The initial guess $X^{0}$ is set to $10^{-5}I$ and the RRN is defined as
\begin{align*}
\text{RRN}:
=\frac{\|X^{k}-X_{\ast}^{0}\|_{F}^{2}}{\|X_{\ast}^{0}\|_{F}^{2}}
=\frac{\|X^{k}-(A^{+}CB^{+}+X^{0}-A^{+}AX^{0}BB^{+})\|_{F}^{2}}
{\|A^{+}CB^{+}+X^{0}-A^{+}AX^{0}BB^{+}\|_{F}^{2}}.
\end{align*}

For these six set of coefficient matrices, the relationship between the RRN and the IT for ME-BK method are given in Fig. 1.
In addition, Table 2 displays some numerical results for the IT and CPU.
It is concluded that the BK method efficiently converges to $X_{\ast}^{0}=A^{+}CB^{+}+X^{0}-A^{+}AX^{0}BB^{+}$ when the matric equation (1.1) is consistent.

\begin{table}[h]
\footnotesize\setlength{\abovecaptionskip}{-0.1cm} 
		\setlength{\belowcaptionskip}{-0.2cm}
\caption{Six sets of coefficient matrices in Example 5.1.}
  \begin{center}
    \begin{tabular}{|c|c|c|c|c|c|c|} 
\hline
\multirow{2}{*}{Matrices}
  &\multicolumn{2}{c|}{Set 1}&\multicolumn{2}{c|}{Set 2}& \multicolumn{2}{c|}{Set 3}\\
\cline{2-7}
 &$A$&$B$&$A$ & $B$&$A$ & $B$\\
\hline
name &\makecell[c]{A=bibd$\_$11$\_$5;\\$\rm A=A^{T}$}&bibd$\_$12$\_$4&\makecell[c]{randn\\(110, 45)}
&\makecell[c]{randn\\(50, 140)}
&\makecell[c]{bibd$\_$12$\_$4} &\makecell[c]{ash958}\\
\hline
size &$462\times 55$&$66\times 495$&$110\times 45$&$50\times 140$&$66\times 495$&$958\times 292$\\
\hline
rank &55&66&45&50&66&292\\
\hline
density &$18.18\%$&$9.09\%$&$100.00\%$&$100.00\%$&$9.09\%$&$0.68\%$\\
\hline
\multirow{2}{*}{Matrices}
  &\multicolumn{2}{c|}{Set 4}&\multicolumn{2}{c|}{Set 5}& \multicolumn{2}{c|}{Set 6}\\
\cline{2-7}
 &$A$&$B$&$A$ & $B$&$A$ & $B$\\
\hline
name&\makecell[c]{randn\\(40, 110)}&\makecell[c]{randn\\(130, 50)}
&\makecell[c]{n3c6-b2}&\makecell[c]{cis-n4c6-b1}
&\makecell[c]{A=randn\\(30, 135);\\A=[A;A]}&\makecell[c]{B=randn\\(120, 50);\\B=[B,B]}\\
\hline
size&$40\times 110$&$130\times 50$
&$455\times 105$&$210\times 21$
&$60\times 135$&$120\times 100$\\
\hline
rank&40&50&91&20&30&50\\
\hline
density&$100.00\%$&$100.00\%$&$2.86\%$&$9.52\%$&$100.00\%$&$100.00\%$\\
\hline
    \end{tabular}
  \end{center}
\end{table}

\begin{table}[h]
\footnotesize\setlength{\abovecaptionskip}{-0.1cm} 
		\setlength{\belowcaptionskip}{-0.2cm}
\caption{The results of IT and CPU for ME-BK method for six sets of coefficient matrices in Table 1.}
  \begin{center}
    \begin{tabular}{|c|c|c|c|c|c|c|c|c|c|c|} 
\hline
\multicolumn{1}{|c|}{Matrices} &Set 1&Set 2&Set 3
&Set 4&Set 5&Set 6\\
\hline
IT&9256&10172&4630&12301&5661&7010\\
\hline
CPU&1.0921&0.2696&13.0959&0.8917&0.5776&0.7468\\
\hline
    \end{tabular}
  \end{center}
\end{table}

\begin{figure}[h]
\renewcommand{\figurename}{Fig.}
\center{\includegraphics[scale=0.6]{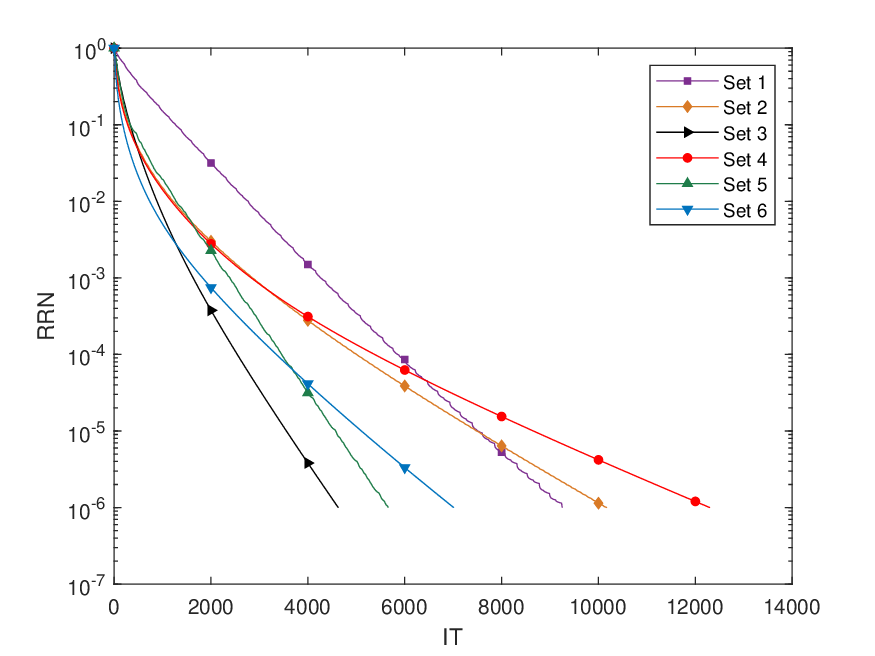}}
\caption{\label{fig:example_figure}{Convergence curves of RRN versus IT for ME-BK method for six set of coefficient matrices in Table 1.}}
\end{figure}

\noindent {\bf Example 5.2.} In this example, we compare the convergence rates of the ME-RBK, ME-GRBK, ME-RGRBK and ME-MWRBK methods. In the following tests, the zero matrix served as an initial guess. The RRN is defined as
\begin{align*}
\text{RRN}:=\frac{\|X^{k}-A^{+}CB^{+}\|_{F}^{2}}{\|A^{+}CB^{+}\|_{F}^{2}},
\end{align*}
and the speed-up of ME-GRBK (ME-RGRBK or ME-MWRBK) method against ME-RBK method is denoted as
\begin{align*}
\text{speed-up}:=\frac{\text{CPU of ME-RBK}}
{\text{CPU of ME-GRBK (ME-RGRBK or ME-MWRBK)}},
\end{align*}
In this example, CPU and IT are the average of 20 independent trials of the relevant method.

In Table 3, we provide six sets of coefficient matrices, where $A$ has full column rank and $B$ has full row rank. In this case, matrix equation (1.1) has a unique solution $(A^{T}A)^{-1}A^{T}CB^{T}(BB^{T})^{-1}$, i.e. $A^{+}CB^{+}=(A^{T}A)^{-1}A^{T}CB^{T}(BB^{T})^{-1}$. Thus, the RRN in Fig. 2 is represented as
\begin{align*}
\rm{RRN}:=\frac{\|X^{k}-(A^{T}A)^{-1}A^{T}CB^{T}(BB^{T})^{-1}\|_{F}^{2}}
{\|(A^{T}A)^{-1}A^{T}CB^{T}(BB^{T})^{-1}\|_{F}^{2}}.
\end{align*}
From Table 4, we observe that the ME-GRBK, ME-RGRBK and ME-MWRBK methods outperform the ME-RBK method in terms of both IT and CPU. The speed-up of ME-GRBK method against ME-RBK method is at least 1.97 (Set 2) and at most 7.49 (Set 6), the speed-up of ME-RGRBK method against ME-RBK method is at least 2.21 (Set 1) and at most attaining 8.39 (Set 6), and the speed-up of ME-MWRBK method against ME-RBK method is at least 2.73 (Set 1) and even attains 13.81 (Set 6).
Fig. 2 shows the convergence curves of the ME-RBK, ME-GRBK, ME-RGRBK and ME-MWRBK methods for the first set of coefficient matrices, where their curves are the outputs of optionally one of the 20 independent experiments of each method. It is intuitively obvious that the number of iterative steps required by the ME-GRBK, ME-RGRBK and ME-MWRBK methods is much less than that needed by the ME-RBK method.

\begin{table}[h]
\footnotesize
\setlength{\abovecaptionskip}{-0.1cm} 
		\setlength{\belowcaptionskip}{-0.2cm}
\caption{Six sets of coefficient matrices with $A$ has full column rank and $B$ has full row rank.}
  \begin{center}
    \begin{tabular}{|c|c|c|c|c|c|c|} 
\hline
\multirow{2}{*}{Matrices}
  &\multicolumn{2}{c|}{Set 1}&\multicolumn{2}{c|}{Set 2}& \multicolumn{2}{c|}{Set 3}\\
\cline{2-7}
 &$A$&$B$&$A$ & $B$&$A$ & $B$\\
\hline
name &ash331& ch5-5-b4 &\makecell[c]{randn\\(140, 30)}&\makecell[c]{randn\\(70, 160)}
&\makecell[c]{sprandn(210,\\ 20, 0.05, 0.1)} &\makecell[c]{sprandn(40,\\ 270, 0.05, 0.1)}\\
\hline
size &$331\times 104$&$120\times 600$&$140\times 30$&$70\times 160$&$210\times 20$&$40\times 270$\\
\hline
rank &104&120&30&70&20&40\\
\hline
density &$1.92\%$&$0.83\%$&$100.00\%$&$100.00\%$&$5.00\%$&$5.00\%$\\
\hline
\multirow{2}{*}{Matrices}
  &\multicolumn{2}{c|}{Set 4}&\multicolumn{2}{c|}{Set 5}& \multicolumn{2}{c|}{Set 6}\\
\cline{2-7}
 &$A$&$B$&$A$ & $B$&$A$ & $B$\\
\hline
name&ash608&bibd$\_$15$\_$3
&\makecell[c]{randn\\(230, 50)}&\makecell[c]{randn\\(110, 240)}
&\makecell[c]{sprandn(160,\\ 35, 0.05, 0.1)} &\makecell[c]{sprandn(25,\\ 150, 0.08, 0.1)}\\
\hline
size&$608\times 188$&$105\times 455$
&$230\times 50$&$110\times 240$
&$160\times 35$&$25\times 150$\\
\hline
rank&188&105&50&110&35&25\\
\hline
density&$1.06\%$&$2.86\%$&$100.00\%$&$100.00\%$&$5.00\%$&$8.00\%$\\
\hline
    \end{tabular}
  \end{center}
\end{table}

\begin{table}[h]
\footnotesize
\small\setlength{\abovecaptionskip}{-0.1cm} 
		\setlength{\belowcaptionskip}{-0.2cm}
\caption{Comparative results of IT and CPU for ME-RBK, ME-GRBK, ME-RGRBK and ME-MWRBK methods for six sets of coefficient matrices in Table 3.}
  \begin{center}
    \begin{tabular}{|c|c|c|c|c|c|c|c|} 
\hline
\multicolumn{2}{|c|}{Matrices} &Set 1&Set 2&Set 3
&Set 4&Set 5&Set 6\\
\hline
\multirow{2}{*}{ME-RBK}& IT&1993.0&9672.6&88971.5&14185.8&14101.1&121604.1\\
\cline{2-8}
 & CPU&1.8595&2.8502&33.1615&32.4754&9.8407&38.6424\\
\hline
\multirow{3}{*}{ME-GRBK}&IT&629.6&4905.5&9597.0&6008.4&8740.3&15452.7\\
\cline{2-8}
 & CPU & 0.9350&1.4482&4.6696&13.5610&3.8892&5.1589\\
 \cline{2-8}
  & speed-up & 1.99&1.97&7.10&2.39&2.53&7.49\\
\hline
\multirow{3}{*}{ME-RGRBK}&IT&579.8&4886.9&9595.0&5997.4&8731.9&15451.0\\
\cline{2-8}
 & CPU&0.8402&1.2686&4.1769&12.2799&3.5242&4.6074\\
 \cline{2-8}
  & speed-up&2.21&2.25&7.94&2.64&2.79&8.39\\
\hline
\multirow{3}{*}{ME-MWRBK}&IT& 542.0&4878.0&9591.5&5996.0&8726.0&15451.0\\
\cline{2-8}
 & CPU & 0.6820&0.9395&2.6801&11.7690&2.8708&2.7980\\
 \cline{2-8}
  & speed-up &2.73&3.03&12.37&2.76&3.43&13.81\\
\hline
    \end{tabular}
  \end{center}
\end{table}

Table 5 shows six sets of coefficient matrices with $A$ has full row rank and $B$ has full column rank. In this case, matrix equation (1.1) has infinitely many solutions and the unique least-norm solution is $A^{T}(AA^{T})^{-1}C(B^{T}B)^{-1}B^{T}$, that is, $A^{+}CB^{+}=A^{T}(AA^{T})^{-1}C(B^{T}B)^{-1}B^{T}$.

\begin{figure}[H]
\renewcommand{\figurename}{Fig.}
\center{\includegraphics[scale=0.6]{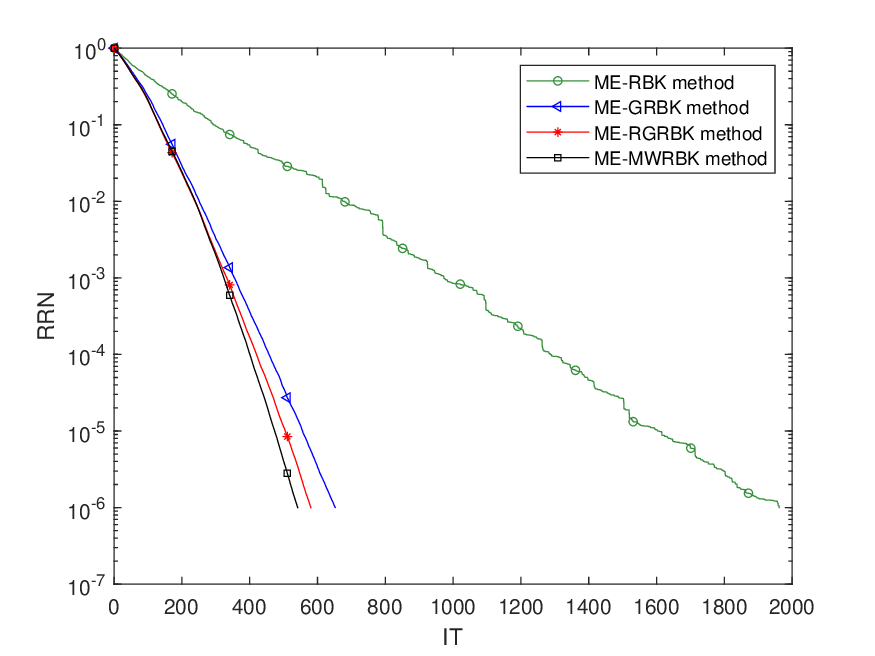}}
\caption{\label{fig:example_figure}{Convergence curves of RRN versus IT for ME-RBK, ME-GRBK, ME-RGRBK and ME-MWRBK methods for the 1-st set of coefficient matrices in Table 3.}}
\end{figure}

\begin{table}[H]
\footnotesize
\small\setlength{\abovecaptionskip}{-0.1cm} 
		\setlength{\belowcaptionskip}{-0.2cm}
\caption{Six sets of coefficient matrices with $A$ has full row rank and $B$ has full column rank.}
  \begin{center}
    \begin{tabular}{|c|c|c|c|c|c|c|} 
\hline
\multirow{2}{*}{Matrices}
  &\multicolumn{2}{c|}{Set 1}&\multicolumn{2}{c|}{Set 2}& \multicolumn{2}{c|}{Set 3}\\
\cline{2-7}
 &$A$&$B$&$A$ & $B$&$A$ & $B$\\
\hline
name &bibd$\_$12$\_$4&ash219&\makecell[c]{randn\\(45, 210)}&\makecell[c]{randn\\(205, 80)}
&\makecell[c]{sprandn(25,\\205,0.08,0.6)} &\makecell[c]{sprandn(125,\\40,0.05,0.5)}\\
\hline
size &$66\times 495$&$219\times 85$&$45\times 210$&$205\times 80$&$25\times 205$&$125\times 40$\\
\hline
rank &66&85&45&80&25&40\\
\hline
density &$9.09\%$&$2.35\%$&$100.00\%$&$100.00\%$&$8.00\%$&$5.00\%$\\
\hline
\multirow{2}{*}{Matrices}
  &\multicolumn{2}{c|}{Set 4}&\multicolumn{2}{c|}{Set 5}& \multicolumn{2}{c|}{Set 6}\\
\cline{2-7}
 &$A$&$B$&$A$ & $B$&$A$ & $B$\\
\hline
name&gams10am&ash331
&\makecell[c]{randn\\(55, 190)}&\makecell[c]{randn\\(165, 75)}
&\makecell[c]{sprandn(130,\\320,0.03,0.5)} &\makecell[c]{sprandn(435,\\170,0.02,0.4)}\\
\hline
size&$114\times 171$&$331\times 104$
&$230\times 50$&$110\times 240$
&$130\times 320$&$435\times 170$\\
\hline
rank&114&104&50&110&130&170\\
\hline
density&$2.09\%$&$1.92\%$&$100.00\%$&$100.00\%$&$5.00\%$&$2.00\%$\\
\hline
    \end{tabular}
  \end{center}
\end{table}

In Table 6, we give the comparative results of IT and CPU for ME-RBK, ME-GRBK, ME-RGRBK and ME-MWRBK methods for six sets of coefficient matrices in Table 5.
It is found that the number of iterative steps and computational time used in the ME-GRBK, ME-RGRBK and ME-MWRBK methods are less than those used in the ME-RBK method.
The speed-up of ME-GRBK method against ME-RBK method is at least 1.04 (Set 4) and at most 2.05 (Set 2), the speed-up of ME-RGRBK method against ME-RBK method is at least 1.04 (Set 4) and at most attaining 2.09 (Set 2), and the speed-up of ME-MWRBK method against ME-RBK method is at least 1.15 (Set 4) and even attains 2.55 (Set 2); see Table 6.
In addition, the convergence curves of the ME-RBK, ME-GRBK, ME-RGRBK and ME-MWRBK methods for the third set of coefficient matrices in Table 5 are displayed in Fig. 3, and their curves are the results of optionally one of the 20 independent tests of each.
The RRN in Fig. 3 is
\begin{align*}
\text{RRN}:=\frac{\|X^{k}-A^{T}(AA^{T})^{-1}C(B^{T}B)^{-1}B^{T}\|_{F}^{2}}
{\|A^{T}(AA^{T})^{-1}C(B^{T}B)^{-1}B^{T}\|_{F}^{2}}.
\end{align*}
We see that the ME-GRBK, ME-RGRBK and ME-MWRBK methods converge more efficiently than the ME-RBK method.

\begin{table}[h]
\footnotesize\setlength{\abovecaptionskip}{-0.1cm} 
		\setlength{\belowcaptionskip}{-0.2cm}
\caption{Comparative results of IT and CPU for ME-RBK, ME-GRBK, ME-RGRBK and ME-MWRBK methods for six sets of coefficient matrices in Table 5.}
  \begin{center}
    \begin{tabular}{|c|c|c|c|c|c|c|c|} 
\hline
\multicolumn{2}{|c|}{Matrices} &Set 1&Set 2&Set 3
&Set 4&Set 5&Set 6\\
\hline
\multirow{2}{*}{ME-RBK}& IT&5090.8&8947.3&718.2&13916.0&19586.8&6035.9\\
\cline{2-8}
 & CPU&3.9997&5.1523&0.1854&8.4072&9.7131&10.1368\\
\hline
\multirow{3}{*}{ME-GRBK}&IT&4569.1&7691.0&524.5&12563.2&15947.7&4554.0\\
\cline{2-8}
 & CPU &3.5971&2.5162&0.1254&8.0620&5.1513&7.0607\\
 \cline{2-8}
  & speed-up & 1.11&2.05&1.48&1.04&1.89&1.44\\
\hline
\multirow{3}{*}{ME-RGRBK}&IT&4568.4&7690.5&519.9&12561.4&15947.0&4553.0\\
\cline{2-8}
 & CPU&3.5001&2.4608&0.1214&8.0513&4.9294&6.9784\\
 \cline{2-8}
  & speed-up&1.14&2.09&1.53&1.04&1.97&1.45\\
\hline
\multirow{3}{*}{ME-MWRBK}&IT&4568.0&7689.0&514.0&12560.0&15946.0&4552.5\\
\cline{2-8}
 & CPU&3.2007&2.0171&0.0933&7.3001&3.9863&6.6089\\
 \cline{2-8}
  & speed-up &1.25&2.55&1.99&1.15&2.44&1.53\\
\hline
    \end{tabular}
  \end{center}
\end{table}

In Table 7, six sets of coefficient matrices with $A$ and $B$ are rank-deficient are given. In this case, matrix equation (1.1) has infinitely many solutions and the unique least-norm solution is $A^{+}CB^{+}$.

Table 8 shows the comparative results of IT and CPU for ME-RBK, ME-GRBK, ME-RGRBK and ME-MWRBK methods for six sets of coefficient matrices in Table 7. We note that the ME-GRBK, ME-RGRBK and ME-MWRBK methods show significant improvements in terms of IT and CPU compared to the ME-RBK method.
As can be seen from Table 8, the speed-up of ME-GRBK method against ME-RBK method is at least 1.35 (Set 4) and at most 2.98 (Set 2), the speed-up of ME-RGRBK method against ME-RBK method is at least 1.38 (Set 4) and at most attaining 3.21 (Set 2), and the speed-up of ME-MWRBK method against ME-RBK method is at least 1.70 (Set 4) and even attains 4.29 (Set 2).
The observations of the second set of coefficient matrices used for testing are intuitively shown in Fig. 4, where we depict the convergence curves of RRN versus IT for the ME-RBK, ME-GRBK, ME-RGRBK and ME-MWRBK methods, with the curves representing the results of optionally one of the 20 independent experiments of each. It is clear that the ME-GRBK, ME-RGRBK and ME-MWRBK methods converge faster than the ME-RBK method as the number of iteration steps increases.

\begin{figure}[h]
\renewcommand{\figurename}{Fig.}
\center{\includegraphics[scale=0.6]{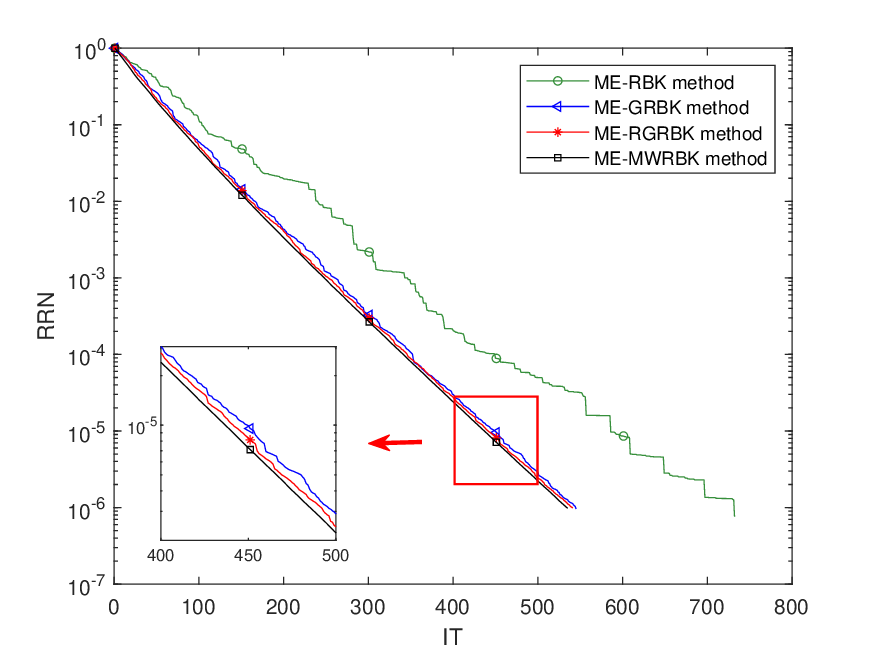}}
\caption{\label{fig:example_figure}{Convergence curves of RRN versus IT for ME-RBK, ME-GRBK, ME-RGRBK and ME-MWRBK methods for the 3-rd set of coefficient matrices in Table 5.}}
\end{figure}

\begin{table}[h]
\footnotesize\setlength{\abovecaptionskip}{-0.1cm} 
		\setlength{\belowcaptionskip}{-0.2cm}
\caption{Six sets of coefficient matrices with $A$ and $B$ are rank-deficient.}
  \begin{center}
    \begin{tabular}{|c|c|c|c|c|c|c|} 
\hline
\multirow{2}{*}{Matrices}
  &\multicolumn{2}{c|}{Set 1}&\multicolumn{2}{c|}{Set 2}& \multicolumn{2}{c|}{Set 3}\\
\cline{2-7}
 &$A$&$B$&$A$ & $B$&$A$ & $B$\\
\hline
name &flower$\_$4$\_$1&n3c6-b2&\makecell[c]{A=randn\\(275, 25);\\ A=[A,A]}&\makecell[c]{B=randn\\(25, 355);\\B=[B;B]}
&\makecell[c]{A=sprandn\\(35,120,0.04,\\0.6); A=[A;A]} &\makecell[c]{B=sprandn\\(115,20,0.03,\\0.5); B=[B,B]}\\
\hline
size &$121\times 129$&$455\times 105$&$275\times 50$&$50\times 355$&$70\times 120$&$115\times 40$\\
\hline
rank &108&91&25&25&35&20\\
\hline
density &$2.47\%$&$2.86\%$&$100.00\%$&$100.00\%$&$4.00\%$&$3.00\%$\\
\hline
\multirow{2}{*}{Matrices}
  &\multicolumn{2}{c|}{Set 4}&\multicolumn{2}{c|}{Set 5}& \multicolumn{2}{c|}{Set 6}\\
\cline{2-7}
 &$A$&$B$&$A$ & $B$&$A$ & $B$\\
\hline
name&n3c6-b1&flower$\_$5$\_$1
&\makecell[c]{A=randn\\(30, 115);\\A=[A;A]}&\makecell[c]{B=randn\\(120, 50);\\B=[B,B]}
&\makecell[c]{A=sprandn\\(205,60,0.05,\\0.5); A=[A,A]} &\makecell[c]{B=sprandn\\(45,220,0.03,\\0.5); B=[B;B]}\\
\hline
size&$105\times 105$&$211\times 201$
&$60\times 115$&$120\times 100$
&$205\times 120$&$90\times 220$\\
\hline
rank&14&179&30&50&60&40\\
\hline
density&$1.90\%$&$1.42\%$&$100.00\%$&$100.00\%$&$5.00\%$&$3.00\%$\\
\hline
    \end{tabular}
  \end{center}
\end{table}

\begin{table}[h]
\footnotesize\setlength{\abovecaptionskip}{-0.1cm} 
		\setlength{\belowcaptionskip}{-0.2cm}
\caption{Comparative results of IT and CPU for ME-RBK, ME-GRBK, ME-RGRBK and ME-MWRBK methods for six sets of coefficient matrices in Table 7.}
  \begin{center}
    \begin{tabular}{|c|c|c|c|c|c|c|c|} 
\hline
\multicolumn{2}{|c|}{Matrices} &Set 1&Set 2&Set 3
&Set 4&Set 5&Set 6\\
\hline
\multirow{2}{*}{ME-RBK}& IT&13588.9&584.8&1043.6&12050.8&5912.0&1847.8\\
\cline{2-8}
 & CPU&8.2546&0.4666&0.2924&4.9501&1.7517&1.2986\\
\hline
\multirow{3}{*}{ME-GRBK}&IT&4828.4&285.5&726.0&7795.5&4865.8&1235.8\\
\cline{2-8}
 & CPU &3.4917&0.1565&0.1532&3.6671&1.2892&0.7054\\
 \cline{2-8}
  & speed-up &2.36&2.98&1.91&1.35&1.36&1.84\\
\hline
\multirow{3}{*}{ME-RGRBK}&IT&4824.2&282.0&724.5&7792.0&4864.6&1235.6\\
\cline{2-8}
 & CPU&3.4075&0.1452&0.1481&3.5757&1.2076&0.7022\\
 \cline{2-8}
  & speed-up&2.42&3.21&1.97&1.38&1.45&1.85\\
\hline
\multirow{3}{*}{ME-MWRBK}&IT&4807.0&272.0&723.0&7778.0&4863.0&1234.0\\
\cline{2-8}
 & CPU &2.9147&0.1087&0.1079&2.9022&0.8829&0.5188\\
 \cline{2-8}
  & speed-up &2.83&4.29&2.71&1.70&1.98&2.50\\
\hline
    \end{tabular}
  \end{center}
\end{table}

\begin{figure}[h]
\renewcommand{\figurename}{Fig.}
\center{\includegraphics[scale=0.6]{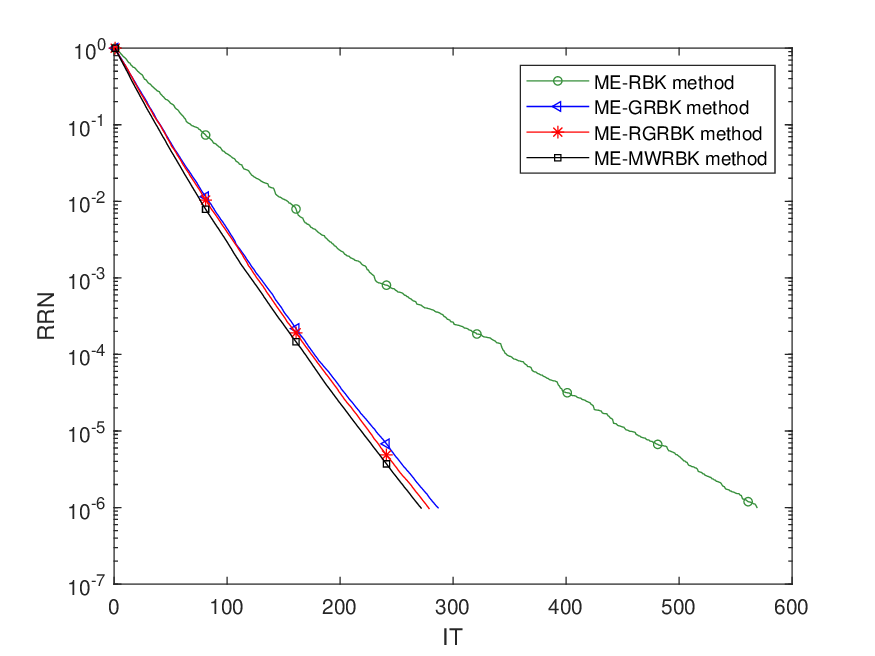}}
\caption{\label{fig:example_figure}{Convergence curves of RRN versus IT for ME-RBK, ME-GRBK, ME-RGRBK and ME-MWRBK methods for the 2-nd set of coefficient matrices in Table 7.}}
\end{figure}

\section{Application in color image restoration}

We first introduce the forward model for color image restoration [5].
The desired original color image $\mathcal{X}$ (the observed color image $\mathcal{B}$) can be shown as a three-dimensional array of size $m\times n\times 3$, where the red (R), green (G) and blue (B) channels are correspondingly three matrices $X_{r}:=\mathcal{X}(:,:,1)$, $X_{g}:=\mathcal{X}(:,:,2)$, $X_{b}:=\mathcal{X}(:,:,3)$ ($B_{r}:=\mathcal{B}(:,:,1)$, $B_{g}:=\mathcal{B}(:,:,2)$, $B_{b}:=\mathcal{B}(:,:,3)$). Let $X=(vec(X_{r}),vec(X_{g}),vec(X_{b}))$ and $B=(vec(B_{r}),vec(B_{g}),vec(B_{b}))$, the forward model is expressed as
\begin{align*}
B=AXA_{c}^{T}+E,
\end{align*}
where $A\in \mathbb{R}^{mn\times mn}$ is the within-channel blurring, $A_{c}\in \mathbb{R}^{3\times 3}$ represents cross-channel blurring and $E\in \mathbb{R}^{mn\times 3}$ is the additive noise. Ideally, if no noise is present in the image generation process, the forward model becomes
\begin{align*}
B=AXA_{c}^{T}.
\end{align*}

Choose ``face" of size $92\times92\times3$, ``bird" of size $96\times96\times3$, ``mandril" of size $125\times120\times3$ and ``barbara" of size $240\times192\times3$
as the test images.
To obtain the blurred images,
the Matlab function `fspecial' was used to generate a Gaussian blur with $G(5,6)$.
Besides, we choose the cross-channel blurring matrix as
\begin{align*}
A_{c}=\left[\begin{array}{cccc}
0.9 &0.05 &0.05\\
0 &0.9 &0.1\\
 0.05 &0.1 &0.85
        \end {array}
\right].
\end{align*}

We use the ME-RBK, ME-GRBK, ME-RGRBK and ME-MWRBK methods in each experiment, and the image restoration results are displayed in Figs. 5-8. The ITs of these experiments in Figs. 5-8 are $5.0\times 10^{4}$, $8.0\times 10^{4}$, $1.0\times 10^{5}$ and $1.5\times 10^{5}$, reapectively.
In all experiments, the parameter $\alpha$ is chosen as $\frac{1}{\|B\|_{2}^{2}}$ and the initial guess is set to a zero matrix.
In the ME-RGRBK method, the relaxation factor $\theta$ is selected to be $0.8$.

To evaluate the restored results, the peak signal-to-noise ratio (PSNR) and the structural similarity index (SSIM) are considered.
The PSNR values of the images obtained after blurring of test images ``face", ``bird", ``mandril" and ``barbara" are 19.3117, 17.5218, 18.8582 and 18.8411 respectively, and their SSIM values are 0.661, 0.660, 0.637 and 0.645 respectively.

\begin{figure}[H]
\renewcommand{\figurename}{Fig.}
\center{\includegraphics[scale=0.9]{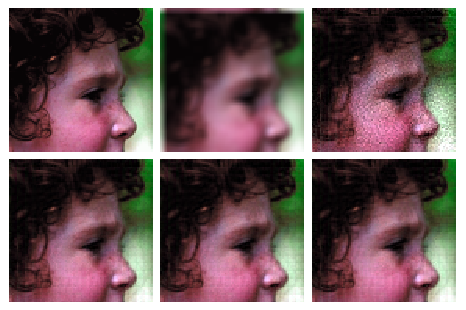}}
\caption{\label{fig:example_figure}{Image restoration results. Top row: test image (left), blurred image (middle), restored image by ME-RBK method (right). Bottom row: restored images by ME-GRBK (left), ME-RGRBK (middle) and ME-MWRBK (right) methods.}}
\end{figure}

\begin{figure}[H]
\renewcommand{\figurename}{Fig.}
\center{\includegraphics[scale=0.9]{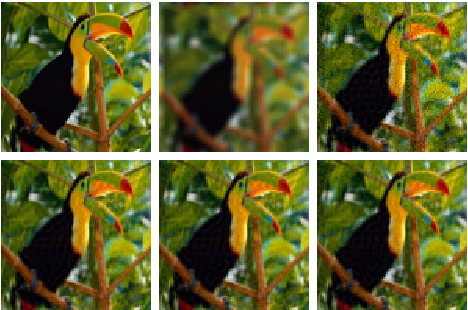}}
\caption{\label{fig:example_figure}{Image restoration results. Top row: test image (left), blurred image (middle), restored image by ME-RBK method (right). Bottom row: restored images by ME-GRBK (left), ME-RGRBK (middle) and ME-MWRBK (right) methods.}}
\end{figure}

\begin{figure}[H]
\renewcommand{\figurename}{Fig.}
\center{\includegraphics[scale=0.8]{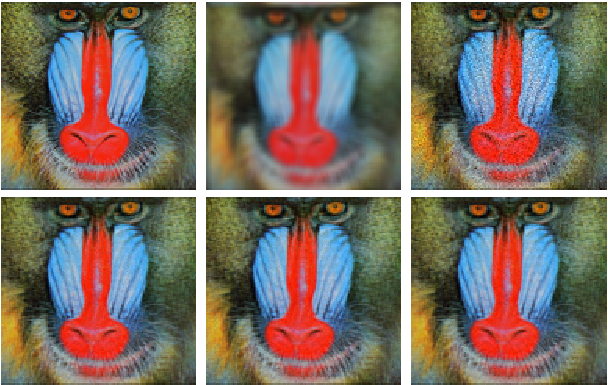}}
\caption{\label{fig:example_figure}{Image restoration results. Top row: test image (left), blurred image (middle), restored image by ME-RBK method (right). Bottom row: restored images by ME-GRBK (left), ME-RGRBK (middle) and ME-MWRBK (right) methods.}}
\end{figure}
\begin{figure}[H]
\renewcommand{\figurename}{Fig.}
\center{\includegraphics[scale=0.5]{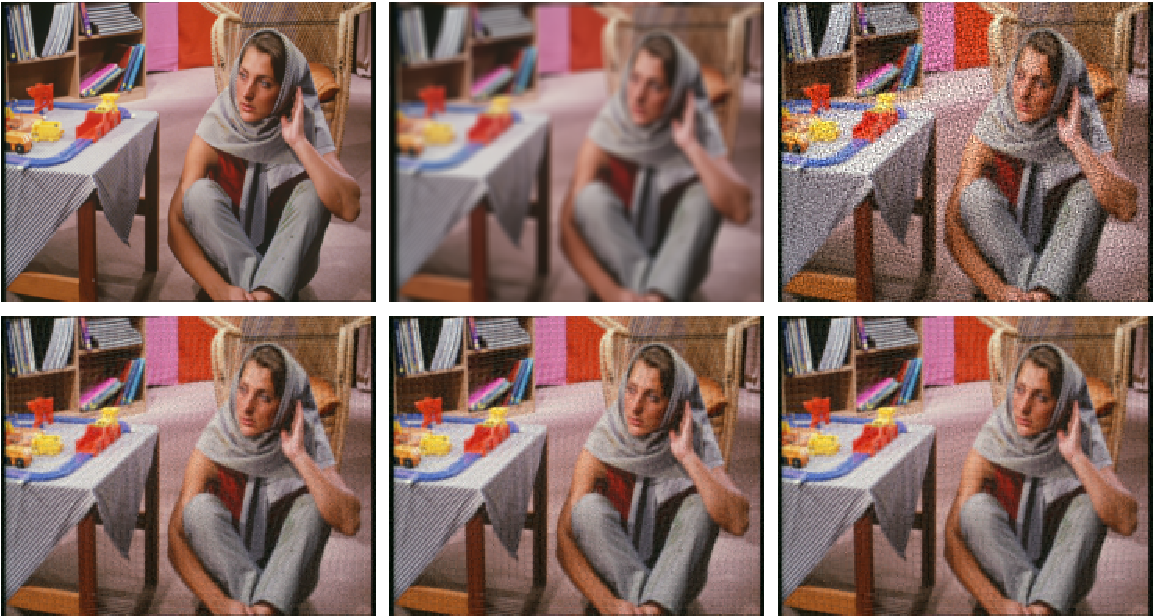}}
\caption{\label{fig:example_figure}{Image restoration results. Top row: test image (left), blurred image (middle), restored image by ME-RBK method (right). Bottom row: restored images by ME-GRBK (left), ME-RGRBK (middle) and ME-MWRBK (right) methods.}}
\end{figure}

Table 9 shows the PSNR and SSIM of the restored images by the ME-RBK, ME-GRBK, ME-RGRBK and ME-MWRBK methods.
The experimental results illustrate that the ME-GRBK, ME-RGRBK and ME-MWRBK methods are more effective than the ME-RBK method in eliminating blurring and obtaining satisfactory restored results.

\begin{table}[H]
\footnotesize\setlength{\abovecaptionskip}{-0.1cm} 
		\setlength{\belowcaptionskip}{-0.2cm}
\caption{The results of PSNR and SSIM of restored images by ME-RBK, ME-GRBK, ME-RGRBK and ME-MWRBK methods.}
  \begin{center}
    \begin{tabular}{|c|c|c|c|c|c|c|c|c|} 
\hline
\multirow{2}{*}{Method}
  &\multicolumn{2}{c|}{face}&\multicolumn{2}{c|}{bird}& \multicolumn{2}{c|}{mandril}&\multicolumn{2}{c|}{barbara}\\
\cline{2-9}
 &PSNR&SSIM&PSNR&SSIM&PSNR&SSIM&PSNR&SSIM\\
\hline
ME-RBK &27.02&0.852&26.85&0.942&25.02&0.869&22.64&0.748\\
\hline
ME-GRBK &33.70&0.939&30.59&0.968&29.90&0.947&30.24&0.940\\
\hline
ME-RGRBK &33.72&0.945&30.61&0.970&29.92&0.948&30.27&0.941\\
\hline
ME-MWRBK &33.72&0.947&30.64&0.971&29.93&0.948&30.29&0.941\\
\hline
    \end{tabular}
  \end{center}
\end{table}

\section{Conclusions}
In this paper, we investigate the ME-GRBK method and its relaxation and deterministic versions, i.e., the ME-RGRBK method and the ME-MWRBK method, to solve the matrix equation (1.1).
It is proved that the ME-GRBK, ME-RGRBK and ME-MWRBK methods  converge to the unique least-norm solution $A^{+}CB^{+}$ when matrix equation (1.1) is consistent.
We also study the ME-BK method and show that it converges to the solution $A^{+}CB^{+}+X^{0}-A^{+}AX^{0}BB^{+}$ when (1.1) is consistent.
Numerical tests demonstrate the effectiveness and superiority of the ME-GRBK, ME-RGRBK and ME-MWRBK methods.
Moreover, the methods in this paper are effectively applied to color image restoration problem, which expands the application range of the Kaczmarz-type method.

\end{document}